\documentclass[10pt]{article}
\usepackage{fullpage,epsfig}
\usepackage{amsmath}
\usepackage{amssymb}
\usepackage{latexsym}
 \numberwithin{equation}{section}
\newtheorem{example}{Example}[section]
\newtheorem{theorem}[example]{Theorem}
\newtheorem{proposition}[example]{Proposition}
\newtheorem{definition}[example]{Definition}
\newtheorem{lemma}[example]{Lemma}

\newtheorem{remark}[example]{Remark}
 
\def\x{{\it x}}

\newcommand{\R}{{\mathbb R}}
\newcommand{\N}{{\mathbb N}}
\newcommand{\ups}{\Upsilon_{{\cal R}}^p(\R^{m\times n})}
\newcommand{\upss}{\Upsilon_{{\cal S}}^p(\R^{m\times n})}

\newcommand{\be}{\begin{eqnarray}}

\newcommand{\ee}{\end{eqnarray}}
\renewcommand{\r}{\varrho}
\renewcommand{\d}{{\rm d}}
\newcommand{\cdm}{{\cal DM}^p_{\cal R}(\O;\R^{m\times n})}
\newcommand{\cdms}{{\cal DM}^p_{\cal S}(\O;\R^{m\times n})}
\newcommand{\gcdm}{{\cal GDM}^p_{\cal R}(\O;\R^{m\times n})}
\newcommand{\dm}{{\cal DM}^p_{\cal R}(\O;\R^{m\times n})}
\newcommand{\gdms}{{\cal GDM}^p_{\cal S}(\O;\R^{m\times n})}

\newcommand{\ra}{\right\rangle}

\newcommand{\la}{\left\langle}
\newcommand{\md}{{\rm d}}

\newcommand{\rems}{{\beta_{{\cal S}}\R^{m\times n}\setminus\R^{m\times n}}}
\renewcommand{\O}{\Omega}
\newcommand{\s}{\sigma}

\renewcommand{\b}{\beta}

\newcommand{\A}[1]{\langle#1\rangle}
\newcommand{\Rn}{\R^{n}}

\newcommand{\wto}{\rightharpoonup}

\newcommand{\B}{{
B}}

\newcommand{\C}{{
C}}

\newcommand{\ITEM}[2]{\parbox[t]{.9cm}{{\rm #1}}\hfill\parbox[t]{153mm}{#2}\vspace*{1mm}\\}

\newcommand{\rca}{{\rm rca}}
\newcommand{\prca}{{\rm rca}_{\vspace*{.5mm}1}^+}

\begin{sloppypar}

\title{Quasiconvexity at the boundary and concentration effects generated by gradients\thanks{This work  was   supported by the grant  P201/10/0357 (GA\v{C}R) }}

\author{ Martin
Kru\v{z}\'{\i}k\thanks{Institute of Information Theory and Automation of the ASCR, Pod vod\'{a}renskou
v\v{e}\v{z}\'{\i}~4, CZ-182~08~Praha~8, Czech Republic (corresponding
address) \& Faculty of Civil Engineering, Czech Technical
University, Th\'{a}kurova 7, CZ-166~ 29~Praha~6, Czech Republic  ({\tt
kruzik@utia.cas.cz})}}

\begin{document}
\date{}
\maketitle

\bigskip

\bigskip

\begin{abstract}
We characterize generalized Young measures, the so-called DiPerna-Majda measures  which are generated by sequences of gradients. In particular,
we precisely describe these measures at the boundary of the domain in the case of the compactification of $\R^{m\times n}$ by the sphere.  We show that this characterization is closely related to the notion of quasiconvexity at the boundary introduced by Ball and Marsden \cite{bama}. As a consequence we get new results on weak $W^{1,2}(\O;\R^3)$ sequential continuity of $u\mapsto a\cdot[{\rm Cof}\, \nabla u]\varrho$, where $\O\subset\R^3$
has a smooth boundary and $a,\varrho$ are certain smooth mappings. 

\end{abstract}

\medskip

\medskip



\section{Introduction}
Oscillations and/or concentrations appear in many problems in the calculus of variations,  partial differential equations, or  optimal
control theory, which admit  only $L^p$ but not $L^\infty$ apriori estimates.
 While Young measures \cite{y} successfully capture oscillatory behavior (see e.g.~\cite{kruzik-luskin,mueller,roubicek-kruzik,roubicek-kruzik-2})  of
sequences they completely miss concentrations. There are several tools how to deal with concentrations.
They can be considered as generalization of Young measures, see for example Alibert's and Bouchitt\'{e}'s approach  \cite{ab}, DiPerna's and Majda's treatment of concentrations \cite{diperna-majda}, or  Fonseca's method described in \cite{fonseca}. An overview can be found in \cite{r,tartar1}.  Moreover, in many cases,
we are interested in oscillation/concentration effects generated by sequences of gradients. A characterization of Young measures generated by gradients was
completely given by Kinderlehrer and Pedregal \cite{k-p1,k-p}, cf. also \cite{mueller, pedregal}.
 The  first attempt to characterize both
oscillations and concentrations in sequences of gradients is due
to Fonseca, M\"{u}ller, and Pedregal \cite{fmp}. They dealt with a special situation of $\{g v(\nabla u_k)\}_{k\in\N}$ where $v$ coincides  with a positively $p$-homogeneous function at infinity (see (\ref{recessionf}) for a precise statement), $u_k\in W^{1,p}(\O;\R^m)$, $p>1$, with  $g$ continuous and   vanishing on $\partial\O$. 
Later on, a characterization of oscillation/concentration effects in terms of  DiPerna's and Majda's generalization
of Young measures was given in \cite{mkak} for arbitrary integrands and in \cite{ifmk} for sequences living in the kernel of a first-order differential operator. Recently Kristensen and Rindler \cite{kristensen-rindler} characterized oscillation/concentration effects in the case $p=1$. 
Nevertheless, a complete  analysis of  boundary effects generated by gradients is still missing. We refer to 
\cite{mkak} for the case where $u_k=u$ on the boundary of the domain. As already observed by Meyers \cite{meyers}, concentration effects at the boundary are closely related to the sequential weak lower semicontinuity of integral functionals $I:W^{1,p}(\O;\R^m)\to\R$: $I(u)=\int_\O v(\nabla u(x))\,\md x$ where  $v:\R^{m\times n}\to\R$ is continuous and such that $|v|\le C(1+|\cdot|^p)$ for some constant $C>0$, cf. also \cite{kroemer} for recent results. Indeed, consider $u\in W^{1,p}_0(B(0,1);\R^m)$, where $B(0,1)$ is the unit ball in $\R^n$ centered at $0$, and extend it by zero to the whole $\R^n$. 
Define for $x\in\R^n$ and $k\in\N$  $u_k(x)=k^{n/p-1}u(kx)$, i.e., $u_k\wto 0$ in $W^{1,p}(B(0,1);\R^m)$ and consider a smooth convex domain $\O\in\R^n$ such that $0\in\partial\O$, $\varrho$ is the outer unit normal to $\partial\O$ at $0$ and let there be $x\in\O$ such that $\varrho\cdot x<0$. Moreover, take a  function  $v$ to be  positively $p$-homogeneous, i.e., $v(\alpha s)=\alpha^p v(s) $ for all $\alpha\ge 0$. Then if $I$ is weakly lower semicontinuous then 
\be\label{example}
0&=&I(0)\le \liminf_{k\to\infty} \int_\O v(\nabla u_k(x))\,\md x= \liminf_{k\to\infty} \int_{B(0,1)\cap\O}v(\nabla u_k(x))\,\md x = \liminf_{k\to\infty}\int_{B(0,1)\cap\O}k^n v(\nabla (u(kx))\,\md x\nonumber \\
&=& \int_{B(0,1)\cap \{x\in\R^n;\ \varrho\cdot x<0\}} v(\nabla u(y))\,\md y\ .
\ee
Thus, we see that $$ 0\le \int_{B(0,1)\cap \{x\in\R^n;\ \varrho\cdot x<0\}} v(\nabla u(y))\,\md y$$ for all $u\in W^{1,p}_0(B(0,1);\R^m)$  forms a necessary condition for weak lower semicontinuity of $I$. 
Here we show that the weak lower semicontinuity of the above defined functional $I$ is intimately related to 
the so-called {\it quasiconvexity at the boundary} defined by Ball and Marsden in \cite{bama} and that this notion of quasiconvexity plays a crucial role in the characterization of parametrized measures generated by sequences of gradients. Moreover, we show that if $\{u_k\}\subset W^{1,2}(\O;\R^3)$, $u_k\wto u$, and $h(x,s):=[{\rm Cof}\, s]\cdot (a(x)\otimes\varrho(x))$ (``Cof'' denotes the cofactor matrix)
for some  $a,\varrho\in C(\bar\O;\R^3)$ such that $\varrho$ coincides with the outer unit normal to $\partial\O$ on the boundary $\partial\O$ of a smooth bounded domain $\O\subset\R^3$ then $h(\cdot,\nabla u_k)\to h(\cdot,\nabla u)$ weakly* in Radon measures supported in $\bar\O$. If, additionally, $h(x,\nabla u_k(x))\ge 0$  for all $k\in\N$ and almost all $x\in\O$ then the above convergence is even in the weak  topology of $L^1(\O)$. Hence, there is a continuous function $\psi:[0,+\infty)\to [0,+\infty)$ such that $\lim_{t\to\infty}\psi(t)/t=+\infty$ and $\sup_{k\in\N}\int_\O\psi\left(h(x,\nabla u_k(x)\right)\,\md x<+\infty$. This result, which can be generalized to higher dimensions, too,  is an analogy  to the  celebrated S.~M\"{u}ler's result on higher integrability of determinants \cite{mueller-det}. See also \cite{Hogan,k-p2}.

\bigskip

\subsection{Basic notation.}
Let us  start with a few definitions and with the explanation of our notation.
Having a bounded domain $\O\subset\R^n$ we denote by $C(\O)$ the space of continuous functions: $\O\to\R$.  Then $C_0(\O)$ consists of  functions from $C(\O)$ whose support is contained in $\O$. In what follows ``{\rm rca}$(S)$'' denotes the set of regular countably additive set functions on the Borel $\s$-algebra on a metrizable set  $S$ (cf. \cite{d-s}), its subset, {\rm rca}$^+_1(S)$,  denotes regular  probability measures on a set $S$.
We write ``$\gamma$-almost all'' or ``$\gamma$-a.e.'' if  we mean ``up to a set with the $\gamma$-measure zero''. If $\gamma$ is the $n$-dimensional Lebesgue measure and $M\subset\R^n$ we omit writing $\gamma$ in the notation.
Further, $W^{1,p}(\O;\R^m)$, $1\le p<+\infty$  denotes the usual space of measurable mappings which are together with
their first (distributional) derivatives integrable with the $p$-th power.
The support of a measure $\sigma\in\ {\rm rca}(\O)$ is a smallest closed set $S$
such that $\sigma(A)=0$ if $S\cap A=\emptyset$. Finally, if $\sigma\in\rca(S)$ we write
$\sigma_s$ and $d_\sigma$ for the singular part and density  of $\sigma$ defined by   the Lebesgue decomposition, respectively. We denote by `w-$\lim$' the weak limit and by $B(x_0,r)$ an open ball in $\R^n$ centered at $x_0$ and the radius $r>0$. The dot product on $\R^n$ is standardly defined as $a\cdot b:=\sum_{i=1}^na_ib_i$ and analogously on $\R^{m\times n}$. Finally, if $a\in\R^m$ and $b\in\R^n$ then  $a\otimes b\in\R^{m\times n}$ with $(a\otimes b)_{ij}=a_ib_j$, and $\mathbb{I}$ denotes the identity matrix.

\bigskip

If not said otherwise, we will  suppose in the sequel  that
$\O\subset\R^n$ is a bounded domain with a $C^1$ boundary. The same regularity is assumed if we say that $\O$ has a smooth boundary.

\bigskip

\subsection{Quasiconvex functions}

Let $\O\subset\R^n$ be a bounded Lipschitz domain. We say that  a function $v:\R^{m\times n}\to\R$ is quasiconvex \cite{morrey} if
for any $s_0\in\R^{m\times n}$ and any $\varphi\in W^{1,\infty}_0(\O;\R^m)$
\be\label{quasiconvexity}
v(s_0)|\O|\le \int_\O v(s_0+\nabla \varphi(x))\,\md x\ .\ee
If $v:\R^{m\times n}\to\R$ is not quasiconvex we define its quasiconvex envelope
$Qv:\R^{m\times n}\to\R$ as
$$
Qv=\sup\left\{h\le v;\ \mbox{$h:\R^{m\times n}\to\R$ quasiconvex }\right\}\ $$
and if the set on the right-hand side is empty we put $Qv=-\infty$.
If $v$ is locally bounded and Borel measurable then for any $s_0\in\R^{m\times n}$ (see \cite{dacorogna})
\be\label{relaxation}
Qv(s_0)=\inf_{\varphi\in W^{1,\infty}_0(\O;\R^m)} \frac{1}{|\O|} \int_\O v(s_0+\nabla \varphi(x))\,\md x\ .\ee

We will also need the following  elementary  result. It can be found in a more general form   e.g. in \cite[Ch.~4, Lemma~2.2]{dacorogna}  or in \cite{morrey}.

\begin{lemma}\label{lemma}  Let $v:\R^{m\times n}\to\R$ be quasiconvex  with  $|v(s)|\le C(1+|s|^p)$, $C>0$, for all $s\in\R^{m\times n}$.
Then there is a constant $\alpha\ge 0$ such that  for every  $s_1,s_2\in\R^{m\times n}$ it holds
\be\label{p-lipschitz-gen}
|v(s_1)-v(s_2)|\le \alpha(1+|s_1|^{p-1}+ |s_2|^{p-1})|s_1-s_2|\ .\ee
\end{lemma}

\bigskip

Following \cite{bama,silhavy, sprenger} we define the notion of  quasiconvexity at the boundary. In order to proceed, we first define the so-called {\it standard boundary domain}.

\begin{definition}\label{def:stabounddom}
Let $\varrho\in\R^n$ be a unit vector and let $\O_\varrho$ be a bounded open Lipschitz domain. We say that 
$\O_\varrho$ is a standard boundary domain with the normal $\varrho$ if there is $a\in\R^n$ such that $\O_\varrho\subset H_{a,\varrho}:=\{x\in\R^n;\ \varrho\cdot x<a\}$ and the $(n-1)$- dimensional interior $\Gamma_\varrho$ of 
$\partial\O_\varrho\cap\partial H_{a,\varrho}$ is nonempty.
\end{definition}

We are now ready to define the quasiconvexity at the boundary.   We put for $1\le p\le+\infty$ 
\be\label{testspace}
W^{1,p}_{\Gamma_\varrho}(\O_\varrho;\R^m):=\{u\in W^{1,p}(\O_\varrho;\R^m);\ u=0 \mbox{ on } \partial\O_\varrho\setminus\Gamma_\varrho\}\ .\ee

\begin{definition} (\cite{bama})\label{qcb-def}
Let $\varrho\in\R^n$ be a unit vector. A function $v:\R^{m\times n}\to\R$ is called quasiconvex at the boundary at $s_0\in\R^{m\times n}$ with respect to $\varrho$ (shortly $v$ is qcb at $(s_0,\varrho)$) if there is $q\in\R^m$ such that for all $u\in W^{1,\infty}_{\Gamma_\varrho}(\O_\varrho;\R^m)$ it holds 
\be\label{qcbinfty}
\int_{\Gamma_\varrho} q\cdot u(x)\,\md S + v(s_0)|\O_\varrho|\le \int_{\O_\varrho} v(s_0+\nabla u(x))\,\md x\ . \ee
\end{definition}
\bigskip

An immediate generalization is the following.

\bigskip

\begin{definition}\label{p-qcb-def}
Let $\varrho\in\R^n$ be a unit vector, $1\le p<+\infty$. A function $v:\R^{m\times n}\to\R$, $|v|\le C(1+|\cdot|^p)$ for some $C>0$  is called $W^{1,p}$-quasiconvex at the boundary at $s_0\in\R^{m\times n}$ with respect to $\varrho$ (shortly $v$ is $p$-qcb at $(s_0,\varrho)$) if there is $q\in\R^m$ such that for all $u\in W^{1,p}_{\Gamma_\varrho}(\O_\varrho;\R^m)$ it holds 
\be\label{qcb}
\int_{\Gamma_\varrho} q\cdot u(x)\,\md S + v(s_0)|\O_\varrho|\le \int_{\O_\varrho} v(s_0+\nabla u(x))\,\md x\ . \ee
\end{definition}

\begin{remark}\label{remark1}
(i) If $v$ is differentiable  then $q:=\frac{\partial v}{\partial s}(s_0)\varrho$ is given uniquely; cf.~\cite{sprenger}.  We denote the set of such of vectors $q$ for which (\ref{qcbinfty}) holds by 
$\partial_v^{\rm qcb}(s_0,\varrho)$. It may be seen as a notion of a ``subdifferential'' for $v$.

\noindent
(ii) It is clear that if $v$ is qcb at $(s_0,\varrho)$ it is also quasiconvex at $s_0$, i.e., (\ref{quasiconvexity}) holds.

\noindent
(iii) If (\ref{qcbinfty}) holds for one standard boundary domain it holds for  other standard boundary domains, too.

\noindent
(iv) If $p>1$,  $v:\R^{m\times n}\to\R$ is positively $p$-homogeneous, i.e.~$v(\lambda s)=\lambda^pv(s)$ for all $s\in\R^{m\times n}$, continuous, and   $p$-qcb at $(0,\varrho)$ then $q=0$ in (\ref{qcbinfty}). Indeed, we have $v(0)=0$ and
suppose that   $\int_{\O_\varrho} v(\nabla u(x))\,\md x<0$ for some $u\in W^{1,\infty}_{\Gamma_\varrho}(\O_\varrho;\R^m)$. 
By (\ref{qcbinfty}), we  must have for all $\lambda>0$
$$ 0\le \lambda^p \int_{\O_\varrho} v(\nabla u(x))\,\md x- \lambda \int_{\Gamma_\varrho} q\cdot u(x)\,\md S\ .$$
However, it is not possible for $\lambda>0$ large enough and therefore for all $u\in W^{1,\infty}_{\Gamma_\varrho}(\O_\varrho;\R^m)$ it holds that 
 $\int_{\O_\varrho} v(\nabla u(x))\,\md x\ge0$. Thus, we can   take  $q=0$.  
\end{remark}

\bigskip

The following lemma shows that Definitions~\ref{qcb-def} and \ref{p-qcb-def} are equivalent for a class of functions whose modulus grows as the $p$-th power.

\begin{lemma} \label{lemma-bdb}
Let $v:\R^{m\times n}\to\R$ be continuous and  such that $|v(A)|\le C(1+|A|^p)$
for all $A\in\R^{m\times n}$ and some  $C>0$ independent of $A$ and some $1\le p<+\infty$. If $v$ is qcb at $(s_0,\varrho)$ it is $p$-qcb at $(s_0,\varrho)$.
\end{lemma}

\bigskip

{\it Proof.}
Take $u\in W_{\Gamma_\varrho}^{1,p}(\O_\varrho;\R^m)$ and a sequence  $\{u_k\}_{k\in\N}\subset W_{\Gamma_\varrho}^{1,\infty}(\O_\varrho;\R^m)$ such that 
$u_k\to u$ strongly in  $W^{1,p}(\O_\varrho;\R^m)$.
We get using (\ref{p-lipschitz-gen}) that 
\begin{align}\label{lim}
\int_{\O_\varrho}v(s_0+\nabla u(x))\,\md x=\lim_{k\to\infty}\int_{\O_\varrho}v(s_0+\nabla u_k(x))\,\md x\ .\end{align}
As $v$ is qcb at $(s_0,\varrho)$  we have
\begin{align}
\lim_{k\to\infty}\int_{\O_\varrho}v(s_0+\nabla u_k(x))\,\md x\ge |\O_{\varrho}|v(s_0)+\lim_{k\to\infty}\int_{\Gamma_\varrho}q\cdot u_k(x)\,\md S
=|\O_{\varrho}|v(s_0)+\int_{\Gamma_\varrho}q\cdot u(x)\,\md S\ ,
\end{align}
which finishes the proof in view of (\ref{lim}).
\hfill
$\Box$
\bigskip

It will be convenient to define the following notion recalling the quasiconvex envelope of $v$ at zero. Here, however, we integrate only over a standard boundary domain with a given normal. If $\varrho\in\R^n$ has a unit length then put
\be\label{envelope}
Q_{b,\varrho}v(0):=\inf_{u\in W_{\Gamma_\varrho}^{1,p}(\O_\varrho;\R^m)}\frac{1}{|\O_\varrho|}\int_{\O_\varrho}v(\nabla u(x))\,\md x\ .\ee

\bigskip

\begin{remark}\label{symmetry}
If $v$ is positively $p$ homogeneous with $p>1$ then  either $Q_{b,\varrho} v(0)=0$ or $Q_{b,\varrho} v(0)=-\infty$.
We also have  that $Q_{b,\varrho}v(0)\le Qv(0)$.  Having a ball $B(0,1)=\{x\in\R^n;\ |x|<1\}$ we put $\O_\varrho:=B(0,1)\cap\{x\in\R^n;\ \varrho\cdot x<0\}$.  In this case, we only integrate over a half-ball in \eqref{envelope}. Hence, we can use only those $u\in  W_{\Gamma_\varrho}^{1,p}(\O_\varrho;\R^m)$ which are symmetric with respect to the plane $\{x;\ \varrho\cdot x=0\}$, i.e., satisfying  $u(x)=u(x-2(\varrho\cdot x)\varrho)$ if $x\in\O_\varrho$.
\end{remark}

\bigskip

Quasiconvexity at the boundary was introduced in \cite{bama} as a necessary condition for strong  local minima of the mixed problem in nonlinear elasticity at boundary points beloging to a free part of the boundary. Contrary to the usual Morrey's quasiconvexity there are not many papers dealing with this notion. Let me point out several interesting results in this direction.  Mielke and Sprenger \cite{mielke-sprenger} investigated relation of quasiconvexity at the boundary and  Agmon's condition for quadratic stored energies in nonlinear elasticity and Sprenger \cite{sprenger} in his thesis defined the so-called polyconvexity at the boundary . Recently, Grabovsky and Mengesha \cite{grabovsky} showed that quasiconvexity at the boundary is a sufficient condition for the so-called $W^{1,\infty}$-sequential weak* local minima - slight weakening of the notion of strong local minimizers. Here we find another interesting connection, namely the fact that quasiconvexity at the boundary plays a crucial role in the analysis of concentration effects generated by gradients and  is essential for  $W^{1,p}$-sequential weak lower semicontinuity of integral functionals as in \eqref{example}.

We start with the following auxiliary lemma. 

\bigskip 

\begin{lemma}\label{lowerb}
Let $v:\R^{m\times n}\to\R$, $|v|\le C(1+|\cdot|^p)$, $C>0$, $1\le p<+\infty$, be quasiconvex, and such that $v(0)= Q_{b,\varrho}v(0)= 0$ for some   $\varrho\in\R^n$.  Let   $\O$ be a bounded domain in $\R^n$ with the $C^1$ boundary, with $x_0\in\partial\O$, and with  $\varrho$ the outer unit normal at $x_0$. Then for every $\varepsilon>0$ there is $\delta>0$   and a continuous function $f:\R\to(0+\infty)$, $\lim_{\varepsilon\to 0}\varepsilon f(\varepsilon)=0$,   such that  $\O\cap B(x_0,\delta)\subset\O$ and it  holds that for every $U\in W^{1,p}_0(\B(0,\delta);\R^m)$ that 
\be
\int_{B(x_0,\delta)\cap\O} v(\nabla U(x))\,\md x\ge   -\varepsilon\int_{B(x_0,\delta)\cap\O}(|\nabla U(x)|+f(\varepsilon)|\nabla U(x)|^p )\,\md x\ ,
\ee 
\end{lemma}

\bigskip

{\it Proof.}
Following \cite{bama}, we can assume, without loss of generality, that $x_0=0$ and that $\varrho=(0,0,\ldots,0, 1)\in\R^n$. Let further 
$$\partial\O\cap B(0,r):=\{x\in B(0,r);\ x_n=h(x')\}\ ,$$
$$\O\cap B(0,r):=\{x\in B(0,r);\ x_n<h(x')\}\ ,$$
where $x=(x_1,\ldots, x_n)\in\R^n$, $x'=(x_1,\ldots, x_{n-1})$, and 
$h\in C^1(\R^{n-1})$ is such that $h(0)=0$ and $\nabla h(0)=0$.
As in \cite{bama} we define for $\xi>0$ $X_\xi:\R^n\to\R^n$ by   $X_\xi(y)=\xi y+h(\xi y')\varrho$. 
 Notice that $\nabla X_\xi(y)=\xi(\mathbb{I}+\varrho\otimes \nabla h(\xi y'))$ and that det$\nabla X_\xi=\xi^n$ because $\varrho\cdot \nabla h=0$. Let $U\in W^{1,p}_0(X_\xi(B(0,r);\R^m)$.
Define $u(y):=\frac1\xi U(X_\xi(y))$, i.e. $u\in W^{1,p}_0(B(0,r);\R^m)$.
Then $\nabla u(y)=\frac1\xi\nabla U(X_\xi(y))\nabla X_\xi(y) = \frac1\xi\nabla U(z)\nabla X_\xi(X^{-1}_\xi(z))$ for $z:=X_\xi(y)$.  Notice that $X^{-1}_{\xi}(z)=\xi^{-1} (z-h( z')\varrho)$ for all $z\in\R^n$. 
 For  $\O_\varrho:=\{\x\in B(0,r);\ x_n<0\}$, we calculate using Lemma~\ref{lemma} 
\begin{align}
\int_{X_\xi(\O_\varrho)}v(\nabla U(z))\,\md z\nonumber\\+\alpha\int_{X_\xi(\O_\varrho)}\left(1+|\nabla U(z)|^{p-1}+\left|\frac1\xi\nabla U(z)\nabla X_\xi(X^{-1}_\xi(z))\right|^{p-1}\right)\left|\nabla U(z)\left(\mathbb{I}- \frac1\xi\nabla X_\xi(X^{-1}_\xi(z)\right)\right|\,\md z\\
\ge \int_{X_\xi(\O_\varrho)} v\left(\frac1\xi\nabla U(z)\nabla X_\xi(X^{-1}_\xi(z))\right)\,\md z=\nonumber\\
=\xi^n\int_{\O_\varrho}v(\nabla u(y))\,\md y\ge 0\ .\nonumber
\end{align}
The last inequality follows from the assumption $Q_{b,\varrho}v(0)\ge 0$. 
Hence, exploiting the identity $\xi^{-1}\nabla X_\xi(X^{-1}_\xi(z))=\mathbb{I}+\varrho\otimes\nabla h(z')$, we get
\begin{align}\label{estimation}
 \int_{X_\xi(\O_\varrho)}v(\nabla U(z))\,\md z \nonumber\\ 
\ge -\alpha\int_{X_\xi(\O_\varrho)}\left(1+|\nabla U(z)|^{p-1}+\left|\frac1\xi\nabla U(z)\nabla X_\xi(X^{-1}_\xi(z))\right|^{p-1}\right)\left|\nabla U(z)\left(\mathbb{I}- \frac1\xi\nabla X_\xi(X^{-1}_\xi(z)\right)\right|\,\md z\\
=-\alpha\int_{X_\xi(\O_\varrho)}(|\nabla U(z)|+|\nabla U(z)|^p(1+|\mathbb{I}+\varrho\otimes\nabla h(z')|^{p-1}))|\varrho\otimes\nabla h(z')|\,\md z \ .\nonumber
\end{align}
However,
\begin{align} 
0\le \int_{X_\xi(\O_\varrho)}(|\nabla U(z)|+|\nabla U(z)|^p(1+|\mathbb{I}+\varrho\otimes\nabla h(z')|^{p-1}))|\varrho\otimes\nabla h(z')|\,\md z\nonumber\\
\le \int_{X_\xi(\O_\varrho)}(|\nabla U(z)|+|\nabla U(z)|^p(1+2^{p-1}(n^{(p-1)/2}+M(|z'|)^{p-1})))M(|z'|)\,\md z\ ,
\end{align}
where $M$ is the modulus of continuity of the uniformly continuous function $z\mapsto\varrho\otimes\nabla h(z')$ on 
$\overline{X_\xi(\O_\varrho)}$, i.e., $\lim_{s\to 0_+}M(s)=0$.
Thus, choosing $\varepsilon>0$, we take $\xi>0$ so small that $\sup_{z\in X_\xi(\O_\varrho)}M(|z|)<\varepsilon/\alpha$ and define $f(\varepsilon):=
1+2^{p-1}(n^{(p-1)/n}+(\varepsilon/\alpha)^{p-1})$.
Then we have 
\begin{align*}0\le\alpha\int_{X_\xi(\O_\varrho)}(|\nabla U(z)|+|\nabla U(z)|^p(1+2^{p-1}(n^{(p-1)/2}+M(|z'|)^{p-1})))M(|z'|)\,\md z\\
\le \int_{X_\xi(\O_\varrho)}(|\nabla U(z)|+f(\varepsilon)|\nabla U(z)|^p )\varepsilon\,\md z\ .
\end{align*}

Finally, in view of (\ref{estimation}) we have 
$$\int_{X_\xi(\O_\varrho)}v(\nabla U(z))\,\md z\ge -\varepsilon\int_{X_\xi(\O_\varrho)}(|\nabla U(z)|+f(\varepsilon)|\nabla U(z)|^p )\,\md z\ .$$
We take $\delta>0$ such that $B(0,\delta)\cap\O\subset X_\xi(\O_\varrho)$ and take $U\in W^{1,p}_0(B(0,\delta);\R^m)$ extended to $\R^n$ by zero which is admissible. Then, as $v(0)=0$, we have from the previous inequality  that 

$$
\int_{B(0,\delta)\cap\O }v(\nabla U(z))\,\md z\ge  -\varepsilon\int_{B(0,\delta)\cap\O}(|\nabla U(z)|+f(\varepsilon)|\nabla U(z)|^p )\,\md z\ .$$

\hfill $\Box$

\bigskip

\begin{example}
If $n=m= 3$  then it is shown in \cite[Prop.~17.2.4]{silhavy} that the function $v:\R^{n\times n}\to\R$ given  by 
$$
v(s)=a\cdot [{\rm Cof} s]\varrho\ $$
is quasiconvex at the boundary with the unit normal $\varrho\in\R^n$. Here $a\in\R^n$ is an arbitrary  constant and ``Cof'' is the cofactor matrix, i.e., ${\rm Cof }s_{ij}=(-1)^{i+j}{\rm det }s'_{ij}$, where $s'_{ij}$ is the submatrix of $s$ obtained from $s$ by removing the $i$-th row and the $j$-th column. Hence,  $v$ is positively $2$-homogeneous. See also \cite{silhavy1} for the role of this $v$ in the definition of 
the so-called interface polyconvexity.  

\end{example}

\bigskip

\subsection{Young measures}
For $p\ge0$ we define the following  subspace of the space
$C(\R^{m\times n})$ of all continuous functions on $\R^{m\times n}$ :
$$
C_p(\R^{m\times n})=\{v\in C(\R^{m\times n}); v(s)=o(|s|^p)\mbox{ for
}|s|\rightarrow\infty\}\ .
$$

The Young
measures on a bounded  domain $\O\subset\Rn$ are weakly* measurable mappings
$x\mapsto\nu_x:\O\to \rca(\R^{m\times n})$ with values in probability measures;
 and the adjective ``weakly* measurable'' means that,
for any $v\in C_0(\R^{m\times n})$, the mapping
$\O\to\R:x\mapsto\A{\nu_x,v}=\int_{\R^{m\times n}} v(\lambda)\nu_x(\d\lambda)$ is
measurable in the usual sense. Let us remind that, by the Riesz theorem,
$\rca(\R^{m\times n})$, normed by the total variation, is a Banach space which is
isometrically isomorphic with $C_0(\R^{m\times n})^*$, where $C_0(\R^{m\times n})$ stands for
the space of all continuous functions $\R^{m\times n}\to\R$ vanishing at infinity.
Let us denote the set of all Young measures by ${\cal Y}(\O;\R^{m\times n})$. It
is known that ${\cal Y}(\O;\R^{m\times n})$ is a convex subset of $L^\infty_{\rm
w}(\O;\rca(\R^{m\times n}))\cong L^1(\O;C_0(\R^{m\times n}))^*$, where the subscript ``w''
indicates the property ``weakly* measurable''.  A classical result
\cite{tartar,warga} is that, for every sequence $\{y_k\}_{k\in\N}$
bounded in $L^\infty(\O;\R^{m\times n})$, there exists its subsequence (denoted by
the same indices for notational simplicity) and a Young measure
$\nu=\{\nu_x\}_{x\in\O}\in{\cal Y}(\O;\R^{m\times n})$ such that
\be\label{jedna2}
\forall v\in C_0(\R^{m\times n}):\ \ \ \ \lim_{k\to\infty}v\circ y_k=v_\nu\ \ \ \
\ \ \mbox{ weakly* in }L^\infty(\O)\ ,
\ee
where $[v\circ y_k](x)=v(y_k(x))$ and
\be
v_\nu(x)=\int_{\R^{m\times n}}v(\lambda)\nu_x(\d\lambda)\ .
\ee
Let us denote by ${\cal Y}^\infty(\O;\R^{m\times n})$ the
set of all Young measures which are created by this way, i.e. by taking
all bounded sequences in $L^\infty(\O;\R^{m\times n})$. Note that (\ref{jedna2}) actually holds for
any $v:\R^{m\times n}\to\R$ continuous.

A generalization of this result was formulated by
Schonbek \cite{schonbek} (cf. also \cite{ball3}): if
$1\le p<+\infty$: for every sequence
$\{y_k\}_{k\in\N}$ bounded in $L^p(\O;\R^{m\times n})$ there exists its
subsequence (denoted by the same
indices) and a Young measure
$\nu=\{\nu_x\}_{x\in\O}\in{\cal Y}(\O;\R^{m\times n})$ such that
\be\label{young}
\forall v\in C_p(\R^{m\times n}):\ \ \ \ \lim_{k\to\infty}v\circ y_k=v_\nu\
\ \ \ \ \ \mbox{ weakly in }L^1(\O)\ .\ee
We say that $\{y_k\}$ generates $\nu$ if (\ref{young}) holds.

  Let us denote by ${\cal
Y}^p(\O;\R^{m\times n})$ the set of all Young measures which are created by this
way, i.e. by taking all bounded sequences in $L^p(\O;\R^{m\times n})$.

The following important lemma was proved in \cite{fmp}.

\begin{lemma}\label{fons}
Let $1 < p<+\infty$ and  $\O\subset\R^n$ be an open bounded set and let $\{ u_k\}_{k\in\N}\subset W^{1,p}(\O;\R^m)$ be bounded. Then there is a subsequence $\{u_j\}_{j\in\N}$ and a sequence $\{z_j\}_{j\in\N}\subset W^{1,p}(\O;\R^m)$ such that
\be\label{rk}
\lim_{j\to\infty} \left|\{x\in\O;\ z_j(x)\ne u_j(x)\mbox{ or }  \nabla z_j(x)\ne \nabla u_j(x)\}\right|=0
\ee
and $\{|\nabla z_j|^p\}_{j\in\N}$ is relatively weakly compact in $L^1(\O)$.
In particular, $\{\nabla u_j\}$ and $\{\nabla z_j\}$ generate the same Young measure.
\end{lemma}

\bigskip

\subsection{DiPerna-Majda measures}
Let us take a complete (i.e. containing constants, separating points
from closed subsets and closed with respect to the Chebyshev norm)
separable  ring ${\cal R}$ of
continuous bounded functions $\R^{m\times n}\to\R$. It is known \cite[Sect.~3.12.21]{engelking} that there is a
one-to-one correspondence ${\cal R}\mapsto\b_{\cal R}\R^{m\times n}$ between such
rings and metrizable compactifications of $\R^{m\times n}$; by a compactification
we mean here a compact set, denoted by $\b_{\cal R}\R^{m\times n}$, into which
$\R^{m\times n}$ is embedded homeomorphically and densely. For simplicity, we
will not distinguish between $\R^{m\times n}$ and its image in $\b_{\cal R}\R^{m\times n}$. Similarly, we will not distinguish between elements of ${\cal R}$ and their unique continuous extensions on $\b_{\cal R}\R^{m\times n}$.

 Let $\s\in\rca(\bar\O)$ be a  positive Radon measure on a bounded domain $\O\subset\R^n$. A
mapping $\hat\nu:x\mapsto \hat\nu_x$ belongs to the
space $L^{\infty}_{\rm w}(\bar{\O},\s;\rca(\b_{\cal R} \R^{m\times n}))$ if it is weakly*  $\s$-measurable (i.e., for any $v_0\in C_0(\R^{m\times n})$, the mapping
$\bar\O\to\R:x\mapsto\int_{\b_{\cal R}\R^{m\times n}} v_0(s)\hat\nu_x(\d s)$ is $\s$-measurable in
the usual sense). If additionally
$\hat\nu_x\in\prca(\b_{\cal R}\R^{m\times n})$ for $\s$-a.a. $x\in\bar\O$
 the collection $\{\hat\nu_x\}_{x\in\bar{\O}}$ is the so-called
Young measure on $(\bar\O,\s)$ \cite{y}, see also
\cite{ball3,r,tartar,valadier,warga}.

DiPerna and Majda \cite{diperna-majda} shown that having a bounded
sequence in $L^p(\O;\R^{m\times n})$ with $1\le p<+\infty$ and
$\O$ an open domain in $\Rn$, there exists its subsequence
(denoted by the same indices), a positive Radon measure
$\s\in\rca(\bar\O)$, and a Young measure  $\hat\nu:x\mapsto
\hat\nu_x$  on  $(\bar\O,\s)$ such that  $(\s,\hat\nu)$ is
attainable by a sequence $\{y_k\}_{k\in\N}\subset
L^p(\O;\R^{m\times n})$ in the sense that $\forall g\!\in\!
C(\bar\O)\ \forall v_0\!\in\!{\cal R}$:
\be\label{basic}\lim_{k\to\infty}\int_\O g(x)v(y_k(x))\d x =
\int_{\bar\O}\int_{\b_{\cal R}\R^{m\times
n}}g(x)v_0(s)\hat\nu_x(\d s)\s(\d x)\ , \ee where \[
v\in\ups:=\{v_0(1+|\cdot|^p);\ v_0\in{\cal R}\}.\]
 In particular,
putting $v_0=1\in{\cal R}$ in (\ref{basic}) we can see that
\be\label{measure} \lim_{k\to\infty}(1+|y_k|^p)\ =\ \s \ \ \ \
\mbox{ weakly* in }\ \rca(\bar\O)\ . \ee If (\ref{basic}) holds,
we say that $\{y_k\}_{\in\N}$ generates $(\sigma,\hat\nu)$. Let us
denote by ${\cal DM}^p_{\cal R}(\O;\R^{m\times n})$ the set of all
pairs $(\s,\hat\nu)\in\rca(\bar\O)\times L^{\infty}_{\rm
w}(\bar{\O},\s; \rca(\b_{\cal R} \R^{m\times n}))$ attainable by
sequences from $L^p(\O;\R^{m\times n})$; note that, taking $v_0=1$
in (\ref{basic}), one can see that these sequences must be
inevitably bounded in $L^p(\O;\R^{m\times n})$. We also denote by $\gcdm$ measures from $\cdm$ generated by a sequence of gradients of some bounded sequence in $W^{1,p}(\O;\R^m)$.  The explicit
description of the elements from ${\cal DM}^p_{\cal
R}(\O;\R^{m\times n})$, called DiPerna-Majda measures, for
unconstrained sequences  was given  in \cite[Theorem~2]{k-r-dm}. 
In fact, it is easy to see that (\ref{basic}) can be also written in the form 
\be\label{basic0}\lim_{k\to\infty}\int_\O h(x,y_k(x))\d x =
\int_{\bar\O}\int_{\b_{\cal R}\R^{m\times
n}}h_0(x,s)\hat\nu_x(\d s)\s(\d x)\ , \ee where $
h(x,s):=h_0(x,s)(1+|s|^p) $ and $h_0\in C(\bar\O\otimes\beta_{\cal R}\R^{m\times n})$. 

We say that $\{y_k\}$ generates $(\sigma,\hat\nu)$ if (\ref{basic}) holds. Moreover, we denote $d_\sigma\in L^1(\O)$ the absolutely continuous (with respect to the Lebesgue measure) part of $\sigma$ in the Lebesgue decomposition of $\sigma$.

Let us recall that for any $(\sigma,\hat\nu)\in {\cal
DM}^p_{\cal R}(\O;\R^{m\times n})$ there is precisely one $(\sigma^\circ,\hat\nu^\circ)\in{\cal
DM}^p_{\cal R}(\O;\R^{m\times n})$ such that
\be\label{nonconc}\int_\O\int_{\R^{m\times n}}v_0(s)\hat\nu_x(\d s)g(x)\s(\d x)=
 \int_{\bar\O}\int_{\R^{m\times n}}
v_0(s)\hat\nu^\circ_x(\d s)g(x)\s^\circ(\d x)\ee
 for any $v_0\in
C_0(\R^{m\times n})$ and any $g\in C(\bar\O)$ and $(\s^\circ,\hat\nu^\circ)$
 is attainable
by a sequence $\{y_k\}_{k\in\N}$ such that the set $\{|y_k|^p;\
k\in\N\}$ is relatively weakly compact in $L^1(\O)$; see
\cite{k-r-dm,r} for details. We call $(\sigma^\circ,\hat\nu^\circ)$
the nonconcentrating modification of $(\sigma,\hat\nu)$.
We call $(\sigma,\hat\nu)\in{\cal
DM}^p_{\cal R}(\O;\R^{m\times n})$ nonconcentrating if
$$\int_{\bar\O}\int_{\beta_{{\cal R}}\R^{m\times n}\setminus\R^{m\times n}}\hat\nu_x(\md s)\sigma(\md x)=0\ .$$
There is a one-to-one correspondence between nonconcentrating DiPerna-Majda measures and Young measures; cf. \cite{r}.

We wish to emphasize  the following fact: if  $\{y_k\}\in L^p(\O;\R^{m\times n})$ generates  $(\sigma,\hat\nu)\in\cdm$ and $\sigma$ is absolutely continuous with respect to the Lebesgue measure it generally {\tt does not } mean that   $\{|y_k|^p\}$ is weakly relatively compact in $L^1(\O)$. A simple examples  can be found e.g. in \cite{k-r-control,r}.

Having a sequence bounded in $L^p(\O;\R^{m\times n})$ generating a
DiPerna-Majda measure $(\sigma,\hat\nu)\in\dm$ it also generates
an $L^p$-Young measure $\nu\in{\cal Y}^p(\O;\R^{m\times n})$.
 It
easily follows from \cite[Th.~3.2.13]{r} that
\be\label{relation}
\nu_x(\md s)=d_{\sigma^\circ}(x)\frac{\hat\nu^\circ_x(\md
s)}{1+|s|^p}\ \  \mbox{ for a.a. $x\in\O$}\ . \ee Note that
(\ref{relation}) is well-defined as $\hat\nu^\circ_x$ is supported
on $\R^{m\times n}$. As pointed out in \cite[Remark~2]{k-r-dm} for
almost all  $x\in\O$ \be\label{density}
d_\sigma(x)=\left(\int_{\R^{m\times n}}\frac{\hat\nu_x(\md
s)}{1+|s|^p}\right)^{-1}\ . \ee In fact, 
that (\ref{nonconc}) can be even improved to
\be\label{nonconc1}\int_\O\int_{\R^{m\times n}}v_0(s)\hat\nu_x(\d
s)g(x)\s(\d x)=
 \int_{\bar\O}\int_{\R^{m\times n}}
v_0(s)\hat\nu^\circ_x(\d s)g(x)\s^\circ(\d x)\ee
 for any $v_0\in {\cal R}$ and any $g\in C(\bar\O)$. The one-to-one correspondence between Young and DiPerna-Majda
 measures, in particular (see (\ref{relation}) and (\ref{nonconc1}))
\[
 \int_{\R^{m\times n}} v(s)\nu_x(ds)=d_\sigma (x) \int_{\R^{m\times n}} v_0(s)\hat{\nu}_x(ds)
 \]
whenever $v\in\ups$,
 finally yields that $\forall g\!\in\! C(\bar\O)\ \forall
v\!\in\ups$: \be\label{basic1}\lim_{k\to\infty}\int_\O
g(x)v(y_k(x))\d x
&=&\int_\O\int_{\R^{m\times n}}v(s) g(x)\nu_x(\md s)\,\md x\nonumber\\
&+&
\int_{\bar\O}\int_{\b_{\cal R}\R^{m\times n}\setminus\R^{m\times n}}\frac{v(s)}{1+|s|^p}\hat\nu_x(\md s)g(x)\s(\md x)\ ,
\ee
where $\nu\in{\cal Y}^p(\O;\R^{m\times n})$ and $(\sigma,\hat\nu)\in\cdm$ are
Young and DiPerna-Majda measures generated by $\{y_k\}_{k\in\N}$, respectively.
We will denote elements from $\cdm$ which are generated by $\{\nabla u_k\}_{k\in\N}$ for some bounded $\{u_k\}\subset W^{1,p}(\O;\R^m)$ by $\gcdm$.

We will also use the following result, whose proof can be found in
several places in various contexts
 (see \cite[Lemma~1., Th.~1,2]{k-r-dm},
 \cite[Prop.~3.2.17]{r}, or for a  compactification of $\R^{m\times n}$ by the sphere  see \cite[Proposition 4.1, part
 (iii)]{ab}).

\begin{lemma}\label{nossin11}
Let $\O\subset\R^n$ be a bounded open domain such that $|\partial\O
| =0$, ${\cal R}$ be a separable complete subring of the ring of
all continuous bounded functions on $\R^{m\times n}$ and
$(\s,\hat{\nu})\in \cdm$. Then for  $\sigma_s$- almost all $x\in\bar\O$
we have \be\label{nossin}
 \hat{\nu}_x(\R^{m\times n})=0.\ee
\end{lemma}

\bigskip

 The following two theorems were proved in \cite{mkak}.

\begin{theorem}\label{suff1}
Let  $\O\subset\R^n$ be a bounded   domain with Lipschitz boundary,
$1<p<+\infty$ and $(\sigma,\hat\nu)\in\cdm$. Then then there is
$u\in W^{1,p}(\O;\R^m)$ and a bounded sequence
$\{u_k-u\}_{k\in\N}\subset W^{1,p}_0(\O;\R^m)$ such that $\{\nabla
u_k\}_{k\in\N}$ generates $(\sigma,\hat\nu)$ if and only if the
following three conditions hold

\be\label{firstmoment6} \mbox{  for a.a. $x\in\O$:  }  \nabla
u(x)=d_\sigma(x)\int_{\beta_{\cal R}\R^{m\times n}}\frac{s}{1+|s|^p}\hat\nu_x(\d
s)\ ,\ee for almost all $x\in\O$ and for all  $v\in\ups$  the
following inequality  is fulfilled \be\label{qc6} Qv(\nabla
u(x))\le d_\sigma(x)\int_{\beta_{\cal R}\R^{m\times
n}}\frac{v(s)}{1+|s|^p}\hat\nu_x(\d s)\ , \ee for $\sigma$-almost
all $x\in\bar\O$ and all $v\in\ups$ with $Qv>-\infty$  it holds
that \be\label{rem6}
 0\le  \int_{\beta_{{\cal R}}\R^{m\times n}\setminus\R^{m\times n}}\frac{v(s)}{1+|s|^p}\hat\nu_x(\md s)\ .
    \ee
\end{theorem}

\bigskip

The next theorem  addresses
DiPerna-Majda measures generated by gradients of maps with
 possibly different  traces.

\bigskip

\begin{theorem}\label{nec50}
Let $\O$ be an arbitrary bounded domain, $1<p<+\infty$ and
$(\sigma,\hat\nu)\in\gcdm$ be  generated by $\{\nabla
u_k\}_{k\in\N}$ such that w-$\lim_{k\to\infty}u_k=u$ in
$W^{1,p}(\O;\R^m)$. Then the conditions (\ref{firstmoment6}),
(\ref{qc6}) hold, and (\ref{rem6}) is satisfied  for $\sigma$-a.a.
$x\in\O$.
\end{theorem}

\bigskip

\begin{remark}\label{superqc}
(i) It can happen that under the assumptions of Theorem~\ref{nec50} formula (\ref{rem6}) does not hold on $\partial\O$.  See an  example
in \cite{murat} showing the violation of weak sequential continuity of $W^{1,2}(\O;\R^2)\to L^1(\O):u\mapsto{\rm det}\ \nabla u$ if $\O=(-1,1)^2$.\\
\noindent
(ii) In terms of Young measures,  conditions (\ref{firstmoment6}) and  (\ref{qc6}) read, respectively:
there is $u\in W^{1,p}(\O;\R^m)$: 
\be\label{firstmoment6ym}
\nabla u(x)=\int_{\R^{m\times n}}s\nu_x(\md s)\ ,\ee 
for all $v:\R^{m\times n}\to\R$, $|v|\le C(1+|\cdot|^p)$:
\be\label{qc6ym}
Qv(\nabla u(x))\le \int_{\R^{m\times n}}v(s)\nu_x(\md s)\ .\ee
\end{remark}

Finally, we have the following result from \cite{mkak}.

 \begin{theorem}\label{wlsc1}
Let $\O\subset\R^n$ be a bounded Lipschitz domain.
Let $0\le g\in C(\bar\O)$, $v\in C(\R^{m\times n})$, $|v|\le
C(1+|\cdot|^p)$, $C>0$, quasiconvex, and  $1<p<+\infty$. Then the
functional $I:W^{1,p}(\O;\R^m)\to\R$ defined as  
\be\label{fcion}I(u):=\int_\O g(x)v(\nabla u(x))\,\md x\ee
 is sequentially weakly
lower semicontinuous in $W^{1,p}(\O;\R^m)$ if and only if for any
bounded sequence $\{w_k\}\subset W^{1,p}(\O;\R^m)$ such that
$\nabla w_k\to 0$ in measure we have
$\liminf_{k\to\infty}I(w_k)\ge I(0)$.
\end{theorem}

\subsubsection{Compactification of $\R^{m\times n}$ by the sphere}

In what follows we will  work   mostly  with a particular compactification of $\R^{m\times n}$, namely, with the compactification by the sphere. We will consider
the following ring  of continuous bounded functions
\begin{eqnarray}\label{spherecomp}
{\cal S}&:=&\left\{^{^{^{^{^{}}}}} v_0\in C(\R^{m\times n}):\mbox{ there exist } c\in\R^{m\times n}\ ,\ v_{0,0}\in C_0(\R^{m\times n}),\mbox{ and }  v_{0,1}\in C(S^{(m\times n)-1}) \mbox{ s.t. }\right.\nonumber\\
&  & \left.
 v_0(s) = c+ v_{0,0}(s)+v_{0,1}\left(\frac{s}{|s|}\right)
\frac{|s|^p}{1+|s|^p}\mbox { if $s\ne 0$ and }  v_0(0)=v_{0,0}(0)\right\}\ ,
\end{eqnarray}  
where $S^{m\times n-1}$ denotes the $(mn-1)$-dimensional unit sphere in $\R^{m\times n}$. Then $\b_{\cal
S}\R^{m\times n}$ is homeomorphic to the unit ball $\overline{B(0,1)}\subset \R^{m\times n}$ via the mapping $d:\R^{m\times n}\to B(0,1)$, $d(s):=s/(1+|s|)$ for all $s\in\R^{m\times n}$. Note that $d(\R^{m\times n})$ is dense in $\overline{B(0,1)}$.

For any $v\in\upss$  there exists a continuous and positively $p$-homogeneous function $v_\infty:\R^{m\times n}\to\R$ (i.e. $v_\infty(\alpha s)=\alpha^p v_\infty(s)$ for all $\alpha\ge 0 $ and $ s\in\R^m$) such that 
\be\label{recessionf}
\lim_{|s|\to\infty}\frac{v(s)-v_\infty(s)}{|s|^p}=0\ .
\ee 

Indeed, if $v_0$ is as in (\ref{spherecomp}) and $v=v_0(1+|\cdot|^p)$  then set
$$v_\infty(s):=\left(c+v_{0,1}\left(\frac{s}{|s|}\right)\right)|s|^p\mbox{ for $s\in\R^{m\times n}\setminus\{0\}$.} $$
By continuity we define $v_\infty(0):=0$. It is easy to see that $v_\infty$ satisfies (\ref{recessionf}).
Such $v_\infty$ is called the {\it recession function} of $v$.
It will be useful to denote
\be\label{vs}
v_\mathcal{S}(s):=(c+v_{0,1})\left(\frac{s}{|s|}\right)\ .
\ee

\bigskip

The following lemma can be found in \cite{fmp, ifmk}.

\begin{lemma}\label{lemma-p}
Let $v\in C(\R^{m\times n})$ be Lipschitz continuous on the unit sphere $S^{m\times n-1}$ and 
$p$-homogeneous, $p\ge 1$. Then  $v$ is $p$-Lipschitz, i.e., there is a constant $\alpha>0$ such that for any $s_1,s_2\in\R^{m\times n}$ it holds 
\be\label{p-lipschitz}
|v(s_1)-v(s_2)|\le \alpha (|s_1|^{p-1}+|s_2|^{p-1})|s_1-s_2|\ .\ee
\end{lemma}

\bigskip

\begin{remark}\label{aboutring}
Notice that ${\cal S}$ contains  all functions  $v_0:=v_{0,0}+ v_\infty/(1+|\cdot|^p)$ where $v_{0,0}\in C_0(\R^{m\times n})$ and $v_\infty:\R^{m\times n}\to\R$ is continuous and positively $p$-homogeneous. 
A  weaker version of Theorem~\ref{nec50} tailored to  the sphere compactification 
was also given in \cite{fmp}.

\end{remark}

If $\{u_k\}\subset W^{1,p}(\O;\R^m)$, $p>1$,  is such that $\{\nabla u_k\}$ generates $(\sigma,\hat\nu)\in\cdms$, 
$\{z_k\}$ is as in Lemma~\ref{fons},  $w_k:=u_k-z_k$ for all $k$, and $v\in C(\R^{m\times n})$ is positively $p$-homogeneous then it follows from 
Lemma~\ref{lemma-p} (and the Stone-Weierstrass theorem on approximations of continuous functions by Lipschitz ones on a compact set)  that for all $g\in C(\bar\O)$

\be\label{concentration-remainder}
\lim_{k\to\infty}\int_\O v(\nabla w_k(x))g(x)\,\md x=\int_{\bar\O}\int_\rems \frac{v(s)}{1+|s|^p}g(x)\hat\nu_x(\md s)\,\sigma(\md x)\ .   
\ee

Indeed, 
let $\{\nabla u_k\}$ generate $(\sigma,\hat\nu)$ and
let $\{z_k\}$ be the sequence constructed in Lemma~\ref{fons}.
Denoting $w_k=u_k-z_k$ for any $k\in\N$ we set $R_k=\{x\in\O;\
\nabla w_k(x)\ne 0\}$. Lemma~\ref{fons} asserts that $|R_k|\to 0$
as $k\to\infty$. We  get from Lemma~\ref{lemma} that  for any
$v\in\ups$
 $p$-homogeneous  and any  $g\in\C(\bar\O)$
\begin{eqnarray}\label{calculation}
& & \left|\int_\O g(x)v(\nabla w_k(x))\,\md x-\int_\O g(x)(v(\nabla u_k(x))-v(\nabla z_k(x)))\,\md x\right| \nonumber\\
&\le&  \|g\|_{C(\bar\O)}\left(\int_{R_k}|v(\nabla u_k(x)-\nabla z_k(x))-v(\nabla u_k(x))|\,\md x+\int_{R_k}|v(\nabla z_k(x))|\,\md x\right)\\
&\le&  C\|g\|_{C(\bar\O)}\int_{R_k}\left[(1+|\nabla u_k(x)-\nabla z_k(x)|^{p-1}+|\nabla u_k|^{p-1})|\nabla z_k(x)|+(1+|\nabla z_k|^p)\right]\,\md x\nonumber\\
&\le &  C'\left(\left(\int_{R_k}|\nabla z_k(x)|^p\,\md
x\right)^{1/p} + \int_{R_k} 1+|\nabla z_k(x)|^p\,\md x+\int_{R_k}
|\nabla z_k(x)|\,\md x\right)\nonumber
\end{eqnarray}
for constants  $C,C'>0$ (which may depend also on ${\rm sup}_k\|
\nabla u_k\|_{L^p(\Omega)}$  and ${\rm sup}_k\| \nabla
z_k\|_{L^p(\Omega)}$ ). The last term goes to zero as $k\to\infty$
because $\{|\nabla z_k|^p\}$ is relatively weakly compact in
$L^1(\O)$ and $|R_k|\to 0$ as $k\to \infty$. This calculation
shows that for $v\in\upss$  we can separate oscillation
and concentration effects of $\{\nabla u_k\}$. Indeed, due to
(\ref{basic1}) we have for any $g\in C(\bar\O)$ and any $v\in\upss$
with $v(0)=0$  that \be\label{hepp}\label{mk} \lim_{k\to\infty}\int_\O v(\nabla
w_k(x))g(x)\,\md x&=&\int_{\bar\O}\int_{\rems}\frac{v(s)}{1+|s|^p}\hat\nu_x(\md s)g(x)\,\sigma(\md x)\\\nonumber
&=&\int_{\bar\O}\int_{\rems}\frac{v_\infty(s)}{1+|s|^p}\hat\nu_x(\md s)g(x)\,\sigma(\md x)\ .\ee 

\bigskip

\section{Main results}

Our main result is the following explicit characterization of DiPerna-Majda measures from $\cdms$ which are generated by gradients.

 \begin{theorem}\label{suff-nec}
Let  $\O\subset\R^n$ be a smooth (at least $C^1$) bounded   domain,
$1<p<+\infty$, and $(\sigma,\hat\nu)\in\cdms$. Then then there is
 a bounded sequence
$\{u_k\}_{k\in\N}\subset W^{1,p}(\O;\R^m)$ such that $\{\nabla
u_k\}_{k\in\N}$ generates $(\sigma,\hat\nu)$ if and only if the
following three conditions hold

\be\label{firstmoment7} \mbox{  for a.a. $x\in\O$:  }  \nabla
u(x)=d_\sigma(x)\int_{\beta_{\cal S}\R^{m\times n}}\frac{s}{1+|s|^p}\hat\nu_x(\d
s)\ ,\ee for almost all $x\in\O$ and for all  $v\in\upss$  the
following inequality  is fulfilled \be\label{qc7} Qv(\nabla
u(x))\le d_\sigma(x)\int_{\beta_{\cal S}\R^{m\times
n}}\frac{v(s)}{1+|s|^p}\hat\nu_x(\d s)\ , \ee for $\sigma$-almost
all $x\in\O$ and all $v\in\upss$ with $Qv_\infty>-\infty$  it holds
that \be\label{rem7}
 0\le  \int_{\beta_{{\cal S}}\R^{m\times n}\setminus\R^{m\times n}}\frac{v(s)}{1+|s|^p}\hat\nu_x(\md s)\ ,
    \ee
and for $\sigma$-almost all $x\in\partial\O$ with the outer unit normal to the boundary  $\varrho(x)$  and all  $v\in\upss$ with $Q_{b,\varrho(x)}v_\infty(0)=0$ it holds that 
\be\label{bd7}
0\le  \int_{\beta_{{\cal S}}\R^{m\times n}\setminus\R^{m\times n}}\frac{v(s)}{1+|s|^p}\hat\nu_x(\md s)\ .
\ee    
\end{theorem}

\bigskip

The following results show that sequential weak lower semicontinuity of $I$ from (\ref{fcion}) puts serious restrictions on $v$.

\bigskip 

\begin{theorem}\label{wlsc2}
Let $\O\subset\R^n$ be a smooth bounded domain and  $1<p<+\infty$.
Let $0\le g\in C(\bar\O)$, $0<g$ on $\partial\O$,  $v\in C(\R^{m\times n})$, and  $|v|\le
C(1+|\cdot|^p)$, $C>0$, quasiconvex such that there is a positively $p$-homogeneous function $v_\infty:\R^{m\times n}\to\R$ satisfying  $\lim_{|s|\to\infty}(v(s)-v_\infty(s))/|s|^p=0$.  Then the
functional $I$ defined by (\ref{fcion}) is sequentially weakly
lower semicontinuous in $W^{1,p}(\O;\R^m)$ if and only if $Q_{b,\varrho}v_\infty(0)=0$ for every $\varrho$ a unit outer normal to $\partial\O$.
\end{theorem}

\bigskip

\begin{theorem}\label{wlsc3}
Let $\O\subset\R^n$ be a smooth bounded domain and  $1<p<+\infty$.
Let $0\le g\in C(\bar\O)$, $0<g$ on $\partial\O$, $v\in C(\R^{m\times n})$, and $|v|\le
C(1+|\cdot|^p)$, $C>0$, quasiconvex such that there is a positively $p$-homogeneous function $v_\infty:\R^{m\times n}\to\R$ satisfying  $\lim_{|s|\to\infty}(v(s)-v_\infty(s))/|s|^p=0$. Let $\{u_k\}\subset W^{1,p}(\O;\R^m)$ weakly converge to $u\in W^{1,p}(\O;\R^m)$. Let $|\nabla u_k|^p\to\sigma$ weakly* in $\rca(\bar\O)$.   Then the
functional $I$ defined by (\ref{fcion}) satisfies $I(u)\le\liminf_{k\to\infty}I(u_k)$ if  $Q_{b,\varrho(x)}v_\infty(0)=0$ for every $\varrho(x)$, a unit outer normal to $\partial\O$ at $x\in\partial\O$, for $\sigma$-a.a.~$x\in\partial\O$. 
\end{theorem}

 \bigskip

\begin{theorem}\label{cof}
Let $\O\subset\R^3$ be a smooth bounded  domain.  Let $\{u_k\}\subset W^{1,2}(\O;\R^3)$ be such that 
$u_k\to u$ weakly in $W^{1,2}(\O;\R^3)$. Let $h(x,s)={\rm Cof }\, s\cdot(a(x) \otimes \varrho(x))$, where $a,\varrho\in C(\bar\O;\R^3)$, $\varrho$ coincides at $\partial\O$ with the outer unit normal to $\partial\O$. Then  for all $g\in C(\bar\O)$  
\be\label{weakcont}
\lim_{k\to\infty}\int_{\O} g(x)h(x,\nabla u_k(x))\,\md x= \int_{\O} g(x)h(x,\nabla u(x))\,\md x\ .\ee
If, moreover, for all $k\in\N$ $h(\cdot,\nabla u_k)\ge 0$ almost everywhere in $\O$  then $h(\cdot,\nabla u_k)\to h(\cdot,\nabla u)$ weakly in $L^1(\O)$.

\end{theorem}

\bigskip

\section{Necessary conditions}

In this section, we show that conditions (\ref{firstmoment7})-(\ref{bd7}) are necessary. In fact, only (\ref{bd7}) needs to be proved because
the other conditions are stated in Theorem~\ref{nec50} which appeared in \cite{mkak}.

We start with the  lemma proved in \cite{mkak}.

\begin{lemma}\label{char}
Let  $(\sigma,\hat\nu)\in\dm$ and  an open  domain
$\omega\subseteq\O$ be  such that $\sigma(\partial\omega)=0$. Let
$\{y_k\}_{k\in\N}$ generate $(\sigma,\hat\nu)$ in the sense
(\ref{basic}).
 Then  for all $v_0\in{\cal R}$ and all $g\in C(\bar\O)$
\be \lim_{k\to\infty}\int_\Omega v(y_k)g(x)\chi_\omega(x)\,\md x=
\int_{\Omega}\int_{\beta_{{\cal R}}\R^{m}}v_0(s)\hat\nu_x(\md
s)g(x)\chi_\omega(x)\,\sigma(\md x)\ ,\ee where $\chi_\omega$ is
the characteristic function of $\omega$ in $\O$.

\end{lemma}

\begin{proposition}
Let $p>1$ and let $(\sigma,\hat\nu)\in{\cal DM}^p_{\cal S}(\O;\R^{m\times n})$ be generated by $\{\nabla u_k\}_k$ where $\{u_k\}_k\in W^{1,p}(\O;\R^m)$ is bounded. Then for $\sigma$-almost all $x\in\partial\O$ it holds that for all $v\in\upss$ with $Q_{b,\varrho(x)}v_\infty(0)=0$ 
\be\label{bounda}
0\le \int_\rems\frac{v(s)}{1+|s|^p}\hat\nu_x(\md s)\ .
\ee

\end{proposition}

\bigskip

{\t Proof.}
Let $\{\nabla u_k\}$ generates $(\sigma,\hat\nu)$, $\{u_k\}\subset W^{1,p}(\O;\R^m)$. We decompose $u_k:=z_k+w_k$ by means of Lemma~\ref{fons}. Then $\nabla w_k\to 0$ in measure and $\{\nabla w_k\}$  carries all the concentrations but no oscillations; cf.~\cite{fmp}. In particular, a simple calculation using (\ref{p-lipschitz}) shows that 
$$\lim_{k\to\infty}\int_\O v_\infty(\nabla w_k(x))\,\md x=\int_{\bar\O}\int_\rems\frac{v_\infty(s)}{1+|s|^p}\hat\nu_x(\md s)\sigma(\md x)\ .$$
Take $x_0\in\partial\O$, $\delta>0$ small enough and such that $\sigma(\partial B(x_0,\delta)\cap\bar\O)=0$. As $Q_{b,\varrho(x_0)}v_\infty(0)=0$ we also have that $Q_{b,\varrho(x_0)}( v_\infty +\varepsilon(|\cdot|+f(\varepsilon)|\cdot|^p))(0)\ge 0$ for $\varepsilon,f(\varepsilon)>0$.  Using Lemmas~\ref{lowerb} and \ref{char}, we have   
\begin{align}
0\le\lim_{k\to\infty} \int_{B(x_0,\delta)\cap\O}v_\infty(\nabla w_k(x))+\varepsilon(|\nabla w_k(x)|+f(\varepsilon)|\nabla w_k(x)|^p)\,\md x\nonumber\\ =\int_{\overline{B(x_0,\delta)\cap\O}}\int_\rems\frac{v_\infty(s)+\varepsilon(|s|+f(\varepsilon)|s|^p)}{1+|s|^p}\hat\nu_x(\md s) \sigma(\md x)\ .
\end{align}
Sending $\varepsilon,\delta\to 0$ and using the Lebesgue-Besicovitch  theorem \cite{evans} we get that for $\sigma$-almost all $x_0\in\partial \O$ it holds that 
\be\label{ptwise1}
0\le \int_\rems \frac{v_\infty(s)}{1+|s|^p}\hat\nu_{x_0}(\md s)=\int_\rems \frac{v(s)}{1+|s|^p}\hat\nu_{x_0}(\md s)\ .\ee

We continue similarly as in \cite{fmp}. The previous calculation
yields the existence of a $\sigma$-null set $E_v\subset\partial\O$ such
that
$$0\le\int_{\beta_{{\cal S}}\R^{m\times n}\setminus\R^{m\times n}}\frac{v(s)}{1+|s|^p}\hat\nu_x(\md s)$$ if $x\not\in E_v$ and $\varrho(x)=\varrho(x_0)=:\rho$.
Let $\{v^k_0\}_{k\in\N}$ be a dense subset of ${\cal S}$, so that
$\{v^k\}_{k\in\N}=\{v_0^k(1+|\cdot|^p)\}_{k\in\N}\subset\upss$.  We
define
$$
E=\bigcup_k\bigcup_{\{j\in\N;\ Q_{b,\rho}(v^k_\infty+(1/j)(1+|\cdot|^p)(0)>-\infty
\}} E_{v^k_\infty+(1/j)(1+|\cdot|^p)}\ .$$ Clearly $\sigma(E)=0$. Fix
$x\in(\O\setminus E)$, $v\in\upss$ such that $Q_{b,\rho}v_\infty(0) >-\infty$ and
choose a subsequence (not relabeled) $\{v_0^{k}\}_{k\in\N}$ such
that
$$v_0^k\to v_0 \mbox { in $C(\beta_{{\cal S}}\R^{m\times n})$  and }  \|v_0^k-v_0\|_{C(\beta_{\cal S}\R^{m\times n})}<\frac{1}{j(k)}\ ,$$
where $j(k)\to \infty$ if $k\to\infty$. We have
\begin{eqnarray*}
v^k(s)+\frac{1}{j(k)}(1+|s|^p) &\ge& v^k(s)+
(1+|s|^p)\|v_0^k-v_0\|_{C(\beta_{\cal S}\R^{m\times n})}\\
&\ge& v^k(s)+|v_0^k(s)-v_0(s)|(1+|s|^p)\ge v(s)\ .
\end{eqnarray*}
Thus, $Q_{b,\rho}(v^k_\infty+\frac{1}{j(k)}(1+|s|^p))+>-\infty$, as well, and
because $x\not\in E$  then $x\not\in
E_{v^k_\infty+(1/j(k))(1+|\cdot|^p)}$ and
\begin{eqnarray*}
0&\le&\lim_{k\to\infty}\int_{\beta_{{\cal S}}\R^{m\times
n}\setminus\R^{m\times n}}
\left(v_0^k(s)+\frac{1}{j(k)}\right)\hat\nu_x(\md s)=
\int_{\beta_{{\cal S}}\R^{m\times n}\setminus\R^{m\times n}}v_0(s)\hat\nu_x(\md s)\\
&=& \int_{\beta_{{\cal S}}\R^{m\times n}\setminus\R^{m\times
n}}\frac{v(s)}{1+|s|^p}\hat\nu_x(\md s) .\end{eqnarray*} \hfill
$\Box$

\section{Sufficient conditions}
The goal of this  section is to show that conditions (\ref{firstmoment7})-(\ref{bd7}) are sufficient, as well.
We will use the following lemma from \cite{fmp} concerning Young measures
from ${\cal Y}^p(\O;\R^{m\times n})$ which are generated by sequences of gradients. 

\bigskip

If $B(0,1)$ is the open unit ball in $\R^n$ centered at zero and $\varrho\in\R^n$ a  unit vector  we denote 
$$B_\varrho=\{x\in\R^n;\ x\in B(0,1)\cap \{x\in\R^n;\ \varrho\cdot x<0\}\}\ $$
and $\partial B_\varrho\supset \Gamma_\varrho =\{x\in\partial B_\varrho;\ \varrho\cdot x=0\}$.

We define two sets of measures:

$$
A^\varrho:=\{\mu\in\rca(\rems);\ \mu\ge 0\ ,\  \la\mu; v_0\ra\ge 0\mbox{ for $v_0\in\mathcal{S}$ if } Q_{b,\varrho}v_\infty(0)=0\} $$
and 
$$
H^\varrho:=\{\bar\delta_{\varrho,\nabla u};\ u\in W_0^{1,p}(B(0,1);\R^m)\}\ ,$$
where for all $v\in$ and  positively $p$-homogeneous

$$
\la\bar\delta_{\varrho,\nabla u}, v_0\ra=|B_\varrho|^{-1}\int_{B_\varrho}v_\mathcal{S}\left(\frac{\nabla u(x)}{|\nabla u(x)|}\right)|\nabla u(x)|^p\,\md x \ .$$

As $\rems\cong S^{m\times n-1}$, the unit sphere in $\R^{m\times n}$ centered at zero we consider both $H^\varrho$ and $A^\varrho$  as sets of measures on the unit sphere. Moreover, $H^\varrho\subset A^\varrho$.  Notice, that  by the definition of $v_{\mathcal S}$ we have 
$$\int_{B_\varrho}v_\mathcal{S}\left(\frac{\nabla u(x)}{|\nabla u(x)|}\right)|\nabla u(x)|^p\,\md x =\int_{B_\varrho}v_\infty(\nabla u(x))\,\md x\ .$$
Moreover, in view of Remark~\ref{symmetry} it is sufficient to consider only $u\in W_0^{1,p}(B(0,1);\R^m)\}$ which are symmetric with respect to the plane $\{x\in\R^n;\ \varrho\cdot x=0\}$ in the definition of $H^\varrho$.

\bigskip

\begin{lemma} Let $n\ge 2$. Then the set $H^\varrho$ is convex. \end{lemma}

\bigskip

{\it Proof.}
Consider $u_1, u_2\in W_0^{1,p}(B(0,1);\R^m)\}$. We need to show that for any $0\le\lambda\le 1$
$\lambda\bar\delta_{\varrho,\nabla u_1}+(1-\lambda)\bar\delta_{\varrho,\nabla u_2}\in H^\varrho$. 
Take $x_0\in B(0,1)\cap\{x\in\R^n;\ \varrho\cdot x=0\}$ such that $|x_0|=1/2$.  Define $\tilde u_1(x):=5^{n/p-1}u_1(5x)$ and
$\tilde u_2(x):=5^{n/p-1}u_1(5(x-x_0))$. We see that $\tilde u_1\in W_0^{1,p}(B(0,1/5);\R^m)$ and $\tilde u_1\in W_0^{1,p}(B(x_0,1/5);\R^m)$, and we may extend them by zero to the whole $\R^n$ (we denote the extension again $\tilde u_1$ and $\tilde u_2$) , so that in particular, $\tilde u_1,\tilde u_2\in  W_0^{1,p}(B(0,1);\R^m)$ and they have disjoint supports. Take $u:=\lambda^{1/p} \tilde u_1+(1-\lambda)^{1/p}\tilde u_2(x)$. Then 
\begin{align*}
\int_{B_\varrho} v(\nabla u(x))\,\md x=\lambda\int_{B_\varrho}v(5^{n/p}\nabla u_1(5x))\,\md x +(1-\lambda)\int_{B_\varrho\cap(B(x_0,1/5)}v(5^{n/p}\nabla u_1(5(x-x_0))\,\md x\\
= \lambda\int_{B_\varrho}v(\nabla u_1(y))\,\md y+(1-\lambda)\int_{B_\varrho} v(\nabla u_2(y))\,\md y\ .
\end{align*}
This proves the claim.
\hfill 
$\Box$

 \bigskip
 
 \begin{remark}
 The case $n=1$ is easy because then quasiconvexity at the boundary reduces to convexity and convex functions are 
 bounded from below by an affine function. Hence, (\ref{rem7}) and (\ref{bd7}) always hold. 
 \end{remark}

\bigskip

\begin{proposition}\label{hahn-banach}
The set $A^\varrho$ is the weak* closure of $H^\varrho$. 
\end{proposition}

\bigskip

{\it Proof.} It is a standard application of the Hahn-Banach theorem. Clearly, $H^\varrho\subset A^\varrho$. Take $v_0\in{\cal S}$ such that 
$\la\mu,v_0\ra\ge a$ for some $a\in\R$ and for all $\mu\in H^\varrho$. Then also 
$$\inf_{u\in W_0^{1,p}(B(0,1);\R^m)} \int_{B_\varrho}v_\infty(\nabla u(x))\,\md x\ge a\ $$
and by $p$-homogeneity we have $ \inf_{u\in W_0^{1,p}(B(0,1);\R^m)}  \int_{B_\varrho}v_\infty(\nabla u(x))\,\md x=0$. Therefore $0\ge a$ and 
$Q_{b,\varrho} v_\infty(0)=0$. By the definition of $A^\varrho$ we see that $\la\pi,v\ra\ge 0\ge a$ for all $\pi\in A^\varrho$.
\hfill
$\Box$

\bigskip

The following two sets of measures were defined in \cite{fmp} 

$$
A:=\{\mu\in\rca(\rems);\ \mu\ge 0\ , \la\mu; v_0\ra\ge 0\mbox{ for $v_0\in{\cal S}$ if } Qv_\infty(0)=0\} $$
and 
$$
H:=\{\bar\delta_{\nabla u};\ u\in W_0^{1,p}(B(0,1);\R^m)\}\ ,$$
where for all $v\in$ and  positively $p$-homogeneous

$$
\la\bar\delta_{\nabla u}, v_0\ra=|B(0,1)|^{-1}\int_{B(0,1)}v_\mathcal{S}\left(\frac{\nabla u(x)}{|\nabla u(x)|}\right)|\nabla u(x)|^p\,\md x \ .$$
 We have the following proposition.
 
 \bigskip
 
 \begin{proposition}\label{closure} (See \cite[Proposition~6.1]{fmp}.) 
 The set $A$ is the weak* closure of $H$.
\end{proposition}

\bigskip

The following theorem formulates sufficient conditions for   $(\sigma,\hat\nu)\in\cdms$ to be generated by gradients.

\begin{theorem}\label{suff-nec1}
Let  $\O\subset\R^n$ be a smooth bounded   domain,
$1<p<+\infty$, and $(\sigma,\hat\nu)\in\cdms$. Then then there is
 a bounded sequence
$\{u_k\}_{k\in\N}\subset W^{1,p}(\O;\R^m)$ such that $\{\nabla
u_k\}_{k\in\N}$ generates $(\sigma,\hat\nu)$ if the
following three conditions hold

\be\label{firstmoment8} \mbox{  for a.a. $x\in\O$:  }  \nabla
u(x)=d_\sigma(x)\int_{\beta_{\cal S}\R^{m\times n}}\frac{s}{1+|s|^p}\hat\nu_x(\d
s)\ ,\ee for almost all $x\in\O$ and for all  $v\in\upss$  the
following inequality  is fulfilled \be\label{qc8} Qv(\nabla
u(x))\le d_\sigma(x)\int_{\beta_{\cal S}\R^{m\times
n}}\frac{v(s)}{1+|s|^p}\hat\nu_x(\d s)\ , \ee for $\sigma$-almost
all $x\in\O$ and all $v\in\upss$ with $Qv_\infty>-\infty$  it holds
that \be\label{rem8}
 0\le  \int_{\beta_{{\cal S}}\R^{m\times n}\setminus\R^{m\times n}}\frac{v(s)}{1+|s|^p}\hat\nu_x(\md s)\ ,
    \ee
and for $\sigma$-almost all $x\in\partial\O$ with the outer unit normal to the boundary  $\varrho(x)$  and all  $v\in\upss$ with $Q_{b,\varrho(x)}v_\infty(0)>-\infty$ it holds that 
\be\label{bd8}
0\le  \int_{\beta_{{\cal S}}\R^{m\times n}\setminus\R^{m\times n}}\frac{v(s)}{1+|s|^p}\hat\nu_x(\md s)\ .
\ee    
\end{theorem}

\bigskip

Before we give the proof we just define a restriction of $(\sigma,\hat\nu)\in\cdm$ to a set $\bar\omega\subset\O$. Naturally, this is a couple 
$(\pi,\hat \gamma)\in{\cal DM}^p_{\cal R}(\omega;\R^m)$ such that $\hat\gamma_x=\hat\nu_x$ if $x\in\bar\omega$ and $\pi$ is the restriction of $\sigma$ to $\bar\omega$.

\bigskip

{\it Proof.}
By the assumptions of the theorem the corresponding Young measure $\nu\in\mathcal{Y}^p(\O;\R^{m\times n})$ given by (\ref{relation}) is generated by gradients of mappings from $W^{1,p}(\O;\R^m)$. By Lemma~\ref{fons} we suppose that this  Young measure is generated by $\{\nabla z_k\}_{k\in\N}$ such that $\{|\nabla z_k|^p\}$ is weakly converging in $L^1(\O;\R^m)$ and $\{z_k\}\subset W^{1,p}(\O;\R^m)$ . Then we look for $\{w_k\}_{k\in\N}\subset W^{1,p}(\O;\R^m)$ such that $\{\nabla w_k\}$ generates  given concentrations but no oscillations. Then $\{\nabla z_k+\nabla w_k\}$ generates the whole DiPerna-Majda measure, see (\ref{hepp}).   If $\sigma(\partial\O)=0$ the proof is exactly the same as in \cite[p.~753]{fmp}. By Lemma~\ref{fons} sought $\{w_k\}$ is such that (i) $w_k\to 0$ weakly in $W^{1,p}(\O;\R^m)$ and (ii) $\nabla w_k\to 0$ in measure. Therefore, it is sufficient to find a sequence of gradients whose Young measure is $\{\delta_0\}_{x\in\O}$ and whose DiPerna-Majda satisfies (\ref{rem8}), (\ref{bd8}), and (\ref{firstmoment8}) and (\ref{qc8}) hold with $u=0$. The proof is divided into two steps. The first step deals with the situation that $\sigma$ only concentrates on the boundary. Two cases are considered - a/ the singular part of $\sigma$ is a weighted sum of Dirac masses and - b/ the general case. The second step assumes that $\sigma$ is arbitrary.

\noindent
(i) Suppose first that $\sigma$ concentrates only at the boundary of $\O$.

  Notice that from (\ref{bd8}) it follows that for $\sigma$-almost all $x\in\partial\O$ $\hat\nu_x\in A^{\varrho(x)}$. Hence, there is a bounded  sequence $\{u_k\}\subset W^{1,p}(B(0,1);\R^m)$, each $u_k$ is symmetric with respect to the plane $\{y\in\R^n;\ \varrho(x)\cdot y=0\}$ and
 $$\int_\rems\frac{v_\infty(s)}{1+|s|^p}\hat\nu_{x}(\md s)=\lim_{k\to\infty}|B_{\varrho(x)}|^{-1}\int_{B_{\varrho(x)}} v_\infty(\nabla u_k(y))\,\md y    \ ,$$
 whenever  $Q_{b,\varrho(x)} v_\infty(0)=0$. 
 By  symmetry,
 there is $\hat\mu_x\in\rca(\rems)$ such that for the same $v$ it holds that 
 $$\int_\rems\frac{v_\infty(s)}{1+|s|^p}\hat\mu_{x}(\md s)=\lim_{k\to\infty}|B_{\varrho(x)}|^{-1}\int_{B(0,1)\setminus B_{\varrho(x)}} v_\infty(\nabla u_k(y))\,\md y    \ .$$
 Thus,\be\label{symetrizace}
 \int_\rems v_0(s)\hat\mu_x(\md s)=\int_\rems v_0\left(s\left(\mathbb{I}-2\varrho(x)\otimes\varrho(x)\right)\right)\hat\nu_x(\md s)\ \ee
 for all $v_0\in{\cal S}$. 
 Altogether, there is a bounded open  $O\subset\R^n$, $\O\subset O$ such that  $(\pi,\hat\gamma)\in{\cal DM}^p_{\cal S}(O;\R^n)$ with 
 \be\label{hatgamma}
 \hat\gamma_x:=\begin{cases}
 \frac12\hat\nu_x+\frac12\hat\mu_x\mbox{ if $x\in\partial\O$},\\
 \delta_0 \mbox{ if $x\in \bar O\setminus\partial\O$}, 
 \end{cases}
 \ee 
 and
 
  \be\label{pi}
 \pi:=\begin{cases}
 2\sigma \mbox{ in $\partial\O$},\\
 \mathcal{L}^n \mbox{ in $ \bar O\setminus\partial\O$}.
 \end{cases}
 \ee 
 This means, in particular, that $\hat\gamma_x\in A$ defined in  Lemma~\ref{closure} and   that $\pi$  is the $n$-dimensional  Lebesgue measure in $ \bar O\setminus\partial\O$. Moreover, $\pi$ does not concentrate on $\partial O$ and by Theorem~\ref{suff1}, see also \cite[Th.~1.1]{fmp}, $(\pi,\hat\gamma)$ is generated by gradients $\{\nabla w_k\}$ where $\{w_k\}\subset W^{1,p}(O;\R^m)$  is bounded.

 \noindent
 a/ Assume first, that that the singular part of $\pi$, $\pi_s$ is equal  to $\sum_{i=1}^N a_i\delta x_i$. We know from Proposition~\ref{closure} that if $Qv_\infty(0)=0$ then 
$$\int_\rems\frac{v_\infty(s)}{1+|s|^p}\hat\gamma_{x_i}(\md s)=\lim_{k\to\infty}|B(0,1)|^{-1}\int_{B(0,1)} v_\infty(\nabla u^{i}_k(x))\,\md x    \ $$
for $1\le i\le N$ and $u_k^{i}\subset W^{1,p}_0(B(0,1);\R^m)$ a bounded sequence in $k$. In view of (\ref{hatgamma}), (\ref{symetrizace}), and Proposition~\ref{hahn-banach}
 we see that $u^i_k$
can be taken symmetric with respect to the plane  $\{y;\ \varrho(x_i)\cdot y=0\}$.
Thus, following \cite{fmp}
$$w_k(x):=k^{n/p-1}|B(0,1)|^{-1/p}\sum_{i=1}^N  a_i^{1/p} u^i_k(k(x-x_i))
$$ 
is such that $\{\nabla(z_k+w_k)\}$ generates $(\pi,\hat\gamma)$ and its restriction to $\O$ generates (by symmetry) $(\sigma,\hat\nu)$.

\noindent    
b/ Take
$l\in\N$. There exists a finite partition ${\cal P}_l =
\{O^{l}_{j}\}_{j=1}^{J(l)}$ of $\bar O$ such that
$O^{l}_{j_1}\bigcap O^{}_{j_2} =\emptyset,\ 1\le j_1< j_2\le
J(l)$ and  all $O^{l}_{j}$ are measurable with
diam$(O^{l}_{j}) <1/l$. Besides, we may suppose that, for any
$l\in\N$, the partition ${\cal P}_{l+1}$ is a refinement of ${\cal
P}_l$, int$(O^{l}_{j})\ne \emptyset$ for all $j$ and that $\sigma$-almost all $x\in\partial\O$ belong to int$(O^{l}_{j})$ for some $j$. Let
$\pi_s$ be the singular part of $\pi$. We set
$a_i^l=\pi_s(O^l_i)$, where $\pi_s$ is the singular part of
$\pi$. Let us put
$$
N(l)=\{1\le j\le J(l);\ a_j^l\ne 0\}\ .$$
 
 If $i\in N(l)$ take  $x_i\in\ {\rm int}(O^{l}_{i})$  such that  $x_i\in\partial\O$ if ${\rm int}(O^{l}_{i})\cap\partial\O\ne\emptyset$. I learned the following argument from Stefan Kr\"{o}mer. We define for $x\in\O^{l}_{i}\cap\partial\O$ rotation matrices $R_{x_il}(x)$ such that $\varrho(x)=R_{x_il}(x)\varrho(x_i)$ for every $x\in\partial\O\cap O^l_i$, hence $R_{x_il}(x_i)=\mathbb{I}$. Notice that if  $v$ is quasiconvex at the boundary at zero with the normal $\varrho$ then $s\mapsto v(sR_{x_il})$ is quasiconvex at the boundary at zero with the normal $R_{x_il}\varrho$.  Define  a measure $(\pi^l,\hat\gamma^l)$ by the formula
$\pi^l(\md x)=d_\pi(x)+\sum_{i\in N(l)}a_i^l\delta_{x_i}$
and \be\label{nua1} \hat\gamma_x^l = \left\{ \begin{array}{ll}
           \hat\gamma_x  & \mbox{ if $x\ne x_i$}\\
            \hat\gamma^l_{x_i} & \mbox {if $x=x_i$}\ ,
                            \end{array}
                   \right.
\ee where supp $\hat\gamma^l_{x_i}\subset\rems$ and for any
$v_0\in{\cal S}$ \be\label{newdef} \int_{\beta_{\cal S}\R^{m\times
n}} v_0(s)\hat\gamma^l_{x_i}(\md
s)=\frac{1}{\pi_s(O^l_i)}\int_{O^l_i}\int_{\beta_{\cal
S}\R^{m\times n}} v_0(sR^{-1}_{x_il}(x))\hat\gamma_x(\md s)\,\pi_s(\md x)\ . \ee
Using Lemma \ref{nossin11}
 we
can equivalently rewrite (\ref{newdef}) as
$$
\int_\rems v_0(sR^{-1}_{x_il})\hat\gamma^l_{x_i}(\md
s)=\frac{1}{\pi_s(O^l_i)}\int_{O^l_i}\int_\rems
v_0(sR^{-1}_{x_il}(x))\hat\gamma_x(\md s)\,\pi_s(\md x)\ .
$$
Theorem~\ref{suff1}  implies that  $(\pi^l,\hat\gamma^l)\in{\cal GDM}_{\cal S}^p(O;\R^m)$. Indeed,
the fact that $(\pi^l,\hat{\gamma}^l)\in\cdm$ is checked by using
Proposition~\ref{characterization}. Moreover, an easy verification
shows that (\ref{firstmoment7}),(\ref{qc7}), (\ref{rem7}) are
also satisfied for $(\pi^l,\hat{\gamma}^l)$. 

Let $\{y_k^l\}_{k\in\N}\subset W^{1,p}(O;\R^m)$ be such that
$\{\nabla y_k^l\}_{k\in\N}$ generates $(\pi^l,\hat\gamma^l)$. We have for any $l\in\N$
\be\label{limita}\lim_{k\to\infty}\int_O (1+|\nabla
y_k^l(x)|^p)\,\md x=\pi^l(\bar O)=\pi(\bar O)\ee and for any
$v_0\in{\cal S}$ and any $g\in C(\bar O)$
\begin{align*}
\lim_{l\to\infty}\left|\int_{\bar O}\int_{\beta_{\cal S}\R^{m\times n}} v_0(s)\hat\gamma^l_{x}(\md s)g(x)\, \pi^l(\md x)-\int_{\bar O}\int_{\beta_{\cal S}\R^{m\times n}} v_0(s)\hat\gamma_{x}(\md s)g(x)\, \pi(\md x)\right|\\
 = \lim_{l\to\infty}\left|\sum_{i\in N(l)}g(x_i)\pi_s(O^l_i)\int_\rems v_0(sR^{-1}_{x_il}(x))\hat\gamma^l_{x_i}(\md s)
 -\int_{\bar O}\int_\rems v_0(s)\hat\gamma_{x}(\md s)g(x)\, \pi_s(\md x)\right|\\
=\lim_{l\to\infty}\left|\sum_{i\in N(l)}\left( \int_{O^l_i}\int_\rems v_0(sR^{-1}_{x_il}(x))\hat\gamma_x(\md s)g(x_i)\,\pi_s(\md x)
- \int_{O^l_i}\int_\rems v_0(s)\hat\gamma_{x}(\md s)g(x)\, \pi_s(\md x)\right)\right|\\
\le \lim_{l\to\infty}\sum_{i\in N(l)}\int_{O^l_i}\int_\rems |v_0(sR^{-1}_{x_il}(x))|\hat\gamma_x(\md s)|g(x)-g(x_i)|\, \pi_s(\md x)\\
+ \lim_{l\to\infty}\sum_{i\in N(l)}\int_{O^l_i}\int_\rems |v_0(sR^{-1}_{x_il}(x))- v_0(s)|\hat\gamma_x(\md s)|g(x)|\, \pi_s(\md x)\\
 \le C\pi_s(\bar O)\lim_{l\to\infty}\left(M_g\left(\frac{1}{l}\right)+M_{v_0}\left(\frac{\tilde C}{l}\right)\right) =0\ ,
\end{align*}
where $|v_0|+|g|\le C$, $\tilde C>0$,
 and
 $M_g$ and $M_{v_0}$ are the moduli of continuity of the uniformly continuous $g\in C(\bar O)$ and $v_0\in C(\rems)$. Here we used the fact that $x\mapsto R(x)$ is continuous for $x\in\partial\O$.
Hence, we get for any $v\in\ups$ and  any $g\in C(\bar O)$
$$
\lim_{l\to\infty}\lim_{k\to\infty}\int_O v(\nabla
y_k^l(x))g(x)\,\md x=\int_{\bar O}\int_{\beta_{\cal S}\R^{m\times
n}} {v_0(s)}\hat\gamma_{x}(\md s)g(x)\, \pi(\md x)\ .
$$ 
However, we know by part (i)\,a/  of the proof that 
 $$
\lim_{l\to\infty}\lim_{k\to\infty}\int_\O v(\nabla
y_k^l(x))g(x)\,\md x=\int_{\bar \O}\int_{\beta_{\cal S}\R^{m\times
n}} {v_0(s)}\hat\nu_{x}(\md s)g(x)\, \sigma(\md x)\ .
$$ The proof of this case  is finished by the diagonalization argument.

\noindent
(ii) Assume that $\sigma$ is arbitrary. 
  Take $1>\varepsilon>0$ and take  $\O(\varepsilon)\subset\O$ 
  such that $\sigma(\O\setminus\overline{\O(\varepsilon)})\to 0$ as $\varepsilon\to 1$.
   Moreover, we suppose that $\sigma(\partial\O(\varepsilon))=0$.
By Theorem~\ref{suff1} there is $\{c_k^\varepsilon\}_{k\in \N}\subset W^{1,p}_0(\O(\varepsilon);\R^m)$ so that $\nabla c^\varepsilon_k$
generates the restriction of $(\sigma,\hat\nu)$ to $\overline{\O(\varepsilon)}$. We can thus extend $c^\varepsilon_k$ by zero to the whole
 $\O$ (without changing the notation).

Let us define 

\be\label{hatbeta}
 \hat\beta^\varepsilon_x:=\begin{cases}
 \hat\nu_x\mbox{ if $x\in\overline{\O(\varepsilon)}$},\\
 \delta_0 \mbox{ if $x\in  \O\setminus\overline{\O(\varepsilon)}$}, \\
 \hat\nu_x\mbox{ if $x\in  \partial \O$}.
 \end{cases}
 \ee 
 and
 
  \be\label{zeta}
 \zeta^\varepsilon:=\begin{cases}
 \sigma \mbox{ in $\overline{\O(\varepsilon)}$},\\
 \mathcal{L}^n \mbox{ in $  \O\setminus\overline{\O(\varepsilon)}$},\\
  \sigma \mbox{ if $x\in  \partial \O$}.
 \end{cases}
 \ee 
 
 Then $(\zeta^\varepsilon,\hat\beta^\varepsilon)\in\gdms$ by the construction from (i) and Theorem~\ref{suff1}. Namely, using (i) we construct
 a sequence of gradients $\{\nabla b^\varepsilon_k\}_k$  generating $(\zeta^\varepsilon,\hat\beta^\varepsilon)$ restricted to $\bar\O\setminus\overline{\O(\varepsilon)}$. This sequence does not concentrate on $\partial\overline{\O(\varepsilon)}$ and can be chosen so, that $\{b^\varepsilon_k\}_k\subset W^{1,p}_0(\O(\varepsilon);\R^m)$. Then using Theorem~\ref{suff1} we have that $\{\nabla b^\varepsilon_k+\nabla c^\varepsilon_k\}_k$ generates  $(\zeta^\varepsilon,\hat\beta^\varepsilon)$. Moreover, $\zeta^\varepsilon$ is majorized by $\sigma$ for every $1\ge\varepsilon>0$, therefore there is a uniform bound on $\{\nabla b^\varepsilon_k+\nabla c^\varepsilon_k\}_k$ in $L^p(\O;\R^{m\times n})$
 independently of $\varepsilon$. We can then shift $b^\varepsilon_k+c^\varepsilon_k$ by a constant (dependent on $k$ and $\varepsilon$) so that its average over $\O$ is zero. The Poincar\'{e} inequality then gives us a uniform bound on  $\{b^\varepsilon_k+c^\varepsilon_k\}\subset W^{1,p}(\O;\R^m)$.
 
 Finally notice that for all $v_0\in{\cal S}$, all $g\in C(\bar\O)$ and $\{\varepsilon^l\}_{l\in\N}\subset(0,1)$, $\lim_{l\to\infty}\varepsilon^l=1$, such that $\sigma(\partial\O(\varepsilon^l))=0$ it holds
 \begin{align*}
 \lim_{l\to\infty}\left|\int_{\bar\O}\int_{\beta_{\cal S}\R^{m\times n}}v_0(s)\hat\nu_x(\md s)g(x)\sigma(\md x)-\int_{\bar\O}\int_{\beta_{\cal S}\R^{m\times n}}v_0(s)\hat\beta^{\varepsilon^l}_x(\md s)g(x)\zeta^{\varepsilon^l}(\md x)\right|\\
 =  \lim_{l\to\infty}\left|\int_{\O\setminus\overline{\O(\varepsilon^l)}}\int_{\beta_{\cal S}\R^{m\times n}}v_0(s)\hat\nu_x(\md s)g(x)\sigma(\md x)-\int_{\O\setminus\O(\varepsilon^l)}v_0(0)g(x)\,\md x\right|\\
 \le \lim_{l\to\infty} C(\sigma(\O\setminus\overline{\O(\varepsilon^l)})+\mathcal{L}^n(\O\setminus\O(\varepsilon^l))=0\ ,
\end{align*}
where $C>0$ is a constant depending on $v_0$ and $g$.
Finally, we finish the proof by a diagonalization argument as ${\cal S}$ and $C(\bar\O)$ are separable.
The theorem is proved. 
 \hfill $\Box$

\section{Proofs of Theorems~\ref{wlsc2},~\ref{wlsc3}, and~\ref{cof}}

{\it Proof of Theorem~\ref{wlsc2}.}
Take $\{u_k\}\subset W^{1,p}(\O;\R^m)$ such that $u_k\to u$ weakly.
 Then multiplying the inequalities in  Theorem~\ref{suff-nec} by $g$ as in the theorem and integrating over $\bar\O$, we have for a subsequence realizing $\liminf_{k\to\infty} I(u_k)$ (not relabeled) and generating $(\sigma,\hat\nu)$  that 
\begin{align*}
\lim_{k\to\infty}I(u_k)=\lim_{k\to\infty}\int_\O v(\nabla u_k(x))g(x))\,\md x=\int_{\bar\O}\int_{\beta_{\cal S}\R^{m\times n}}\frac{v(s)}{1+|s|^p}\hat\nu_x(\md s)g(x)\sigma(\md x)
\ge I(u)\ ,\nonumber
 \end{align*}
 which finishes the proof of the ``if part''.
 
To show the ``only if part'' of the assertion we assume that $I$ is sequentially weakly lower semicontinuous and want to show that $Q_{b,\varrho}v_\infty(0)=0$.  Consider $x_0\in\partial\O$ and  $u\in W^{1,p}_0(B(0,1);\R^m)$ and extend it by zero to the whole $\R^n$. 
Define for $x\in\R^n$ and $k\in\N$  $u_k(x)=k^{n/p-1}u(k(x-x_0))$, i.e., $u_k\wto 0$ in $W^{1,p}(\O;\R^m)$ and assume that  $\varrho$ is the outer unit normal to $\partial\O$ at $x_0$. The Young measure generated by $\{\nabla u_k\}$ is just $\{\delta_0\}_{x\in\O}$ because $\nabla u_k\to 0$ in measure.
Hence, we get 
\begin{align}\label{weaklsc}
\lim_{k\to\infty}I(u_k)=I(0)+\int_{\bar\O}\int_{\rems}\frac{v(s)}{1+|s|^p}\hat\nu_x(\md s)g(x)\sigma(\md x)\nonumber\\
=I(0)+\int_{\bar\O}\int_{\rems}\frac{v_\infty(s)}{1+|s|^p}\hat\nu_x(\md s)g(x)\sigma(\md x)\ .
\end{align} 
Defining $I_\infty(u):=\int_\O g(x)v_\infty(\nabla u(x))\,\md x$ we have
\begin{align}
\lim_{k\to\infty}I_\infty(u_k)=\int_{\bar\O}\int_{\rems}\frac{v_\infty(s)}{1+|s|^p}\hat\nu_x(\md s)g(x)\sigma(\md x)\nonumber\\
=g(x_0)\int_{B(0,1)\cap\{x\in\R^n;\,\varrho\cdot x\le 0\}}v_\infty(\nabla u(x))\,\md x\ .
\end{align}
If the last term is negative for some $u\in W_0^{1,p}(B(0,1);\R^m)$, i.e., if  $Q_{b,\varrho}v_\infty(0)=-\infty$ then $\lim_{k\to\infty}I(u_k)<I(0)$, so that 
$I$ is not sequentially weakly lower semicontinuous due to (\ref{weaklsc}) which   contradicts our assumption. The assertion is proved. 
 
\hfill $\Box$

\bigskip

{\it Proof of Theorem~\ref{wlsc3}.} It is just a corollary of Theorem~\ref{wlsc2}.
\hfill $\Box$

\bigskip

{\it Proof of Theorem~\ref{cof}.} Let $m=n=3$. Functions $s\mapsto\pm{\rm Cof}\ s$ are both quasiconvex \cite{dacorogna} and, as already mentioned in  \cite{silhavy}, $s\mapsto \pm a\cdot[{\rm Cof}\ s]\varrho$ is quasiconvex at the boundary with the normal $\varrho$.  Thus, all the inequalities in Theorem~\ref{suff-nec} are equalities if applied to $v=\pm h(x,\cdot)$ for a fixed $x\in\bar\O$. Hence, if $\{\nabla u_k\}$ generates $(\sigma,\hat\nu)$ on the domain $\bar\O$ and $u_k\to u$ weakly in $W^{1,2}(\O;\R^3)$
we have for all $g\in C(\bar\O)$
$$
\lim_{k\to\infty}\int_{\O}g(x)h(x,\nabla u_k(x)) \,\md x=\int_{\bar\O}\int_{\beta_{\cal S}\R^{m\times n}}\frac{h(x,s)}{1+|s|^2}\hat\nu_x(\md s)g(x)\sigma(\md x)= \int_{\O}g(x) h(x,\nabla u(x))\,\md x\ .$$
If $0\le h(x,\nabla u_k(x))$ the result follows by Lemma~\ref{equiintegrability} in the Appendix.
\hfill $\Box$

\bigskip
\appendix
\section{Appendix}

The following proposition from \cite{k-r-dm} explicitly characterizes
elements of $\cdm$.

\bigskip

\begin{proposition}\label{characterization}
Let $\O\subset\R^n$ be a bounded open domain,
${\cal R}$ be a separable complete subring of the ring of all
continuous bounded functions on $\R^{m\times n}$ and $(\s,\hat{\nu})\in
\rca(\bar{\O})\times L^{\infty}_{\rm w}(\bar{\O},\s;
\rca(\b_{\cal R}\R^{m\times n}))$ and $1\le p<+\infty$.
Then the following two statements are equivalent with each other:\\
\ITEM{(i)}{the pair $(\s,\hat\nu)$ is the DiPerna-Majda measure, i.e.
$(\s,\hat\nu)\in{\cal DM}^p_{\cal R}(\O;\R^{m\times n})$,}
\ITEM{(ii)}{The following properties are satisfied simultaneously:
\begin{enumerate}
\item
$\s$ is positive,
\item
$\s_{\hat\nu}\in\rca(\bar\O)$ defined by
$\s_{\hat\nu}(\d x)=(\int_{\R^{m\times n}}\hat\nu_x(\d s))\s(\d x)$
is
absolutely\\
 continuous with respect to the Lebesgue
measure\\
 ($d_{\s_{\hat\nu}}$ will denote its density),
\item for a.a. $x\in\O$ it holds
$$
\!\!\!\!\!\!\!\!\!\!\!\!\!\!\!\!\!\!\!\!\!\!\!\!\!\!\!\!\!\!\!\!\!
\!\!\!\!\!\!\!\!\!\!\!\!\!\!\!\!\!\!\!\!\!\!\!\!\!\ \int_{\R^{m\times n}}\hat\nu_x(\d s) >0,\ \ \ \ \ \ d_{\s_{\hat\nu}}(x)
=\left(\int_{\R^{m\times n}}\frac{\hat\nu_x(\d s)}
{1+|s|^p}\right)^{-1}\int_{\R^{m\times n}}\hat\nu_x(\d s)\ ,
$$
\item
for $\s$-a.a. $x\in\bar\O$ it holds
$$
\hat\nu_x\ge 0,\ \ \ \ \ \
\int_{\b_{\cal R}\R^{m\times n}}\hat\nu_x(\d s)=1\ .
$$
\end{enumerate}
}
\end{proposition}

\begin{lemma}\label{equiintegrability}
Let $1\le p<+\infty$,  $0\le h_0\in C(\bar\O\times\beta_{\cal R}\R^{m\times n})$, and  let  $\{\nabla u_k\}_{k\in\N}\subset L^p(\O;\R^{m\times n})$ with $u_k\in W^{1,p}(\O;\R^m)$,  generate
$(\sigma,\hat\nu)\in {\cal DM}^p_{\cal R}(\O;\R^m)$.  Let $h(x,s):=h_0(x,s)(1+|s|^p)$. Then $\{h(x,\nabla u_k)\}_{k\in\N}$ is weakly relatively compact in $L^1(\O)$ if and only if 
\be\label{wrc}
\int_{\bar\O}\int_{\beta_{\cal R}\R^{m\times n}\setminus\R^{m\times n}} h_0(x,s)\hat\nu_x(\md s)\sigma(\md x)=0\ .
\ee  

\end{lemma}

\bigskip

{\it Proof.}
We follow the proof of \cite[Lemma 3.2.14(i)]{r}. Suppose first that (\ref{wrc}) holds. For  $\r \ge 0$ define  the function $\xi^\r:\R^{m\times n}\to\R$ 
$$
\xi^\r(s): = \left\{ \begin{array}{ll}
                            0 & \mbox{if $|s|\le \r$,}\\
                            |s|-\r & \mbox{if $\r
                                          \le |s |\le \r+1$,}\\
                            1 & \mbox{if $|s |\ge \r+1$}\ .
                            \end{array}
                   \right. 
$$
Note that always $\xi^\r\in{\cal R}$, hence $\xi^\r h_0(x,\cdot)\in {\cal R}$ because ${\cal R}$ is closed under multiplication.   
We have due to the Lebesgue Dominated Convergence Theorem
$$
\lim_{\r\to\infty}\int_{\bar\O}\int_{\b_{\cal R}\R^{m\times n}\setminus B(0,\r)}h_0(x,s)\hat\nu_x(\md s)\sigma(\md x)=\int_{\bar\O}\int_{\beta_{\cal R}\R^{m\times n}\setminus\R^{m\times n}} h_0(x,s)\hat\nu_x(\md s)\sigma(\md x)=0\ .$$
 Let  $\varepsilon>0$ and $\r$ be large enough so that
$$
\int_{\bar\O}\int_{\b_{\cal R}\R^{m\times n}}\xi^\r(s)h_0(x,s)\hat\nu_x(\md s)\sigma(\md x)\le 
\int_{\bar\O}\int_{\b_{\cal R}\R^{m\times n}\setminus B(0,\r)}h_0(x,s)\hat\nu_x(\md s)\sigma(\md x)\le\frac{\varepsilon}{2}\ , $$
and choose $k_\r\in\N$ such that, if $k\ge k_\r$, then 
$$
\left| \int_{\bar\O}\int_{\b_{\cal R}\R^{m\times n}}\xi^\r(s)h_0(x,s)\hat\nu_x(\md s)\sigma(\md x)-\int_\O \xi^\r(\nabla u_k(x))h(x,\nabla u_k(x))\,\md x\right|\le \frac{\varepsilon}{2}\ . $$

Therefore, if $k\ge k_\r$ then 
$\int_\O \xi^\r(\nabla u_k(x))h(x,\nabla u_k(x))\,\md x\le\varepsilon$, and so
$$\int_{\{x\in\O:\ |\nabla u_k(x)|\ge\r+1\}}h(x,\nabla u_k(x))\,\md x\le \int_\O \xi^\r(\nabla u_k(x))h(x,\nabla u_k(x))\,\md x\le\varepsilon\ .$$
As $0\le h(x,\cdot)\le C(1+|\cdot|^p)$ for some $C>0$, we get for $K\ge C(1+(\r+1)^p)$ that 
$$\int_{\{x\in\O:\ |h(x,\nabla u_k(x))|\ge K\}} h(x,\nabla u_k(x))\,\md x\le \int_{\{x\in\O:\ |\nabla u_k(x)|\ge\r+1\}}h(x,\nabla u_k(x))\,\md x\le\varepsilon\ .$$
Clearly, the  finite set $\{h(x,\nabla u_k)\}_{k=1}^{k_\r}$  is weakly relatively compact in $L^1(\O)$, which means that for $K_0>0$ sufficiently large and $1\le k\le k_\r$ 
$$
\int_{\{x\in\O:\ |h(x,\nabla u_k(x))|\ge K_0\}} h(x,\nabla u_k(x))\,\md x\le\varepsilon\ .$$
Hence,  
$$
\sup_{k\in\N}\int_{\{x\in\O:\ |h(x,\nabla u_k(x))|\ge \max(K_0,K)\}} h(x,\nabla u_k(x))\,\md x\le\varepsilon\ ,$$
and $\{h(x,\nabla u_k)\}$ is relatively weakly compact in $L^1(\O)$ by the Dunford-Pettis criterion. 
Consequently, if $\{h(x,\nabla u_k)\}$ is relatively weakly compact in $L^1(\O)$, then the limit of a (sub)sequence can be fully  described by the Young measure generated by $\{\nabla u_k\}$, see e.g.~\cite{ball3,mueller,pedregal}. Hence, $\hat\nu$ is supported on $\R^{m\times n}$.
\hfill
$\Box$

\bigskip

{\bf Acknowledgment:} I thank Agnieszka Ka\l amajska, Stefan Kr\"{o}mer, Alexander Mielke, and Filip Rindler for helpful discussions.

\end{sloppypar}


\bigskip\bigskip

\vspace*{1cm}

\bigskip

\begin{thebibliography}{19}
\baselineskip=12pt
{\footnotesize

\bibitem{ab}
{\sc Alibert, J.J., Bouchitt\'{e}, G.}: Non-uniform integrability and generalized Young measures. {\it. J. Convex Anal.} {\bf 4} (1997), 125--145.


\bibitem{ball3}
{\sc Ball, J.M.}: A version of the fundamental theorem for Young measures. In:
{\it PDEs and Continuum Models of Phase Transition.} (Eds. M.Rascle, D.Serre,
M.Slemrod.) Lecture Notes in Physics {\bf 344}, Springer, Berlin, 1989,
pp.207--215.

\bibitem{bama}
{\sc Ball, J.M., Marsden, J.}: Quasiconvexity at the boundary, positivity of the second variation and elastic stability. {\it Arch.~Rat.~Mech.~Anal.} {\bf 86}(1984), 251--277. 

\bibitem{murat}
{\sc Ball, J.M., Murat, F.}: $W^{1,p}$-quasiconvexity and variational problems for multiple integrals. {\it J. Funct. Anal.} {\bf 58} (1984), 225--253.

\bibitem{zhang}
{\sc Ball, J.M., Zhang K.-W.}: Lower semicontinuity of multiple integrals and the biting lemma. {\it Proc. Roy. Soc. Edinburgh } {\bf 114A} (1990), 367--379.

\bibitem{chacon}
{\sc Brooks, J.K., Chacon, R.V.}: Continuity and compactness in measure.  {\it Adv. in Math.} {\bf 37} (1980), 16--26.
\bibitem{dacorogna}
{\sc Dacorogna, B.} {\it Direct Methods in the Calculus of Variations.}
Springer, Berlin, 1989.

\bibitem{diperna-majda}
{\sc DiPerna, R.J., Majda, A.J.}: Oscillations and concentrations in weak
solutions of the incompressible fluid equations.
{\bf 108} (1987), 667--689.

\bibitem{d-s}
{\sc Dunford, N., Schwartz, J.T.}: {\it Linear Operators.}, Part I,
Interscience, New York, 1967.


\bibitem{engelking}
{\sc Engelking, R.}: {\it General topology} $2^{nd}$ ed., PWN, Warszawa, 1985.

\bibitem{evans}
{\sc Evans, L.C., Gariepy, R.F.}: {\it Measure Theory and Fine Properties of Functions.} CRC Press, Inc. Boca Raton, 1992.
\bibitem{fonseca}
{\sc Fonseca, I.}: Lower semicontinuity of surface energies. {\it Proc. Roy. Soc. Edinburgh} {\bf 120A} (1992), 95--115.

\bibitem{ifmk}
{\sc Fonseca, I., Kru\v{z}\'{\i}k, M.}: Oscillations and concentrations generated by ${\mathcal A}$-free
mappings and weak lower semicontinuity of integral functionals. {\it ESAIM Control Optim.~Calc.~Var.} {\bf 16} (2010), 472--502.

\bibitem{fmp}
{\sc Fonseca, I., M\"{u}ller, S., Pedregal, P.}: Analysis of concentration and oscillation effects generated by gradients. {\it SIAM J. Math. Anal.} {\bf 29} (1998), 736--756.


\bibitem{grabovsky}
{\sc Grabovsky, Y., Mengesha, T.}: Direct approach to the problem of strong local minima in calculus of variations. {\it Calc.~Var.} {\bf 29} (2007), 59--83.

\bibitem{Hogan}
{\sc Hogan, J., Li, C., McIntosh, A., Zhang, K.}: Global higher integrability  of Jacobians on bounded domains.
{\it Annales de l'I.H.P.~section C} {\bf 17} (2000), 193--217.

\bibitem{mkak}
{\sc Ka\l amajska, A., Kru\v{z}\'{\i}k, M.}: Oscillations and concentrations in sequences of gradients.  {\it ESAIM Control Optim.~Calc.~Var.} {\bf 14} (2008),  71-104.


\bibitem{k-p1}
{\sc Kinderlehrer, D., Pedregal, P.}: Characterization of Young measures generated by gradients. {\it Arch. Rat. Mech. Anal.} {\bf 115} (1991), 329--365.

\bibitem{k-p2}
{\sc Kinderlehrer, D., Pedregal, P.}: Weak convergence of integrands and the Young measure representation. {\it SIAM J. Math. Anal.} {\bf 23} (1992), 1--19. 

\bibitem{k-p}
{\sc Kinderlehrer, D., Pedregal, P.}: Gradient Young measures
generated by sequences in Sobolev spaces. {\it J. Geom. Anal.} {\bf4}
(1994), 59--90.



\bibitem{kristensen-rindler}
{\sc Kristensen, J., Rindler, F.}: Characterization of generalized gradient Young measures
generated by sequences in $W^{1,1}$ and BV. {\it  Arch.~Rat.~Mech. Anal.} {\bf 197} (2010), 539--598. 


\bibitem{kroemer}
{\sc Kr\"{o}mer, S.}: On the role of lower bounds in characterizations of weak lower semicontinuity of multiple integrals. {\it Adv.~Calc.~Var.} {\bf 3} (2010), 378--408.

\bibitem{kruzik-luskin}
{\sc Kru\v{z}\'{\i}k, M., Luskin, M.}: The computation of martensitic microstructure with piecewise laminates. {\it J. Sci. Comp.} {\bf 19} (2003), 293--308.


\bibitem{k-r-dm}
{\sc Kru\v{z}\'{\i}k, M., Roub\'{\i}\v{c}ek, T.}:
On the measures of DiPerna and Majda. {\it  Mathematica Bohemica}, {\bf 122} (1997),  383--399.

\bibitem{k-r-control}
{\sc Kru\v{z}\'{\i}k, M., Roub\'{\i}\v{c}ek, T.}: Optimization problems with
concentration and oscillation effects: relaxation theory and numerical
approximation. {\it Numer. Funct. Anal. Optim.}  {\bf 20} (1999), 511-530.

\bibitem{meyers}
{\sc Meyers, N.G.}: Quasi-convexity and lower semicontinuity of
multiple integrals of any order, {\it Trans. Am. Math. Soc.} {\bf
119} (1965), 125-149.

\bibitem{mielke-sprenger}
{\sc Mielke, A., Sprenger, P.}: Quasiconvexity at the boundary and a simple variational formulation of Agmon's condition. {\it J.~Elasticity} {\bf 51} (1998), 23--41.



\bibitem{morrey}
{\sc Morrey, C.B.}: {\it Multiple Integrals in the Calculus of Variations.} Springer, Berlin, 1966.

\bibitem{mueller-det}
{\sc M\"{u}ller, S.}: Higher integrability of determinants and weak convergence in $L^1$. {\it J. reine angew. Math.} {\bf 412} (1990), 20--34.

\bibitem{mueller}
{\sc M\"{u}ller, S.}: {\it Variational models for microstructure and phase transisions.} Lecture Notes in Mathematics {\bf 1713} (1999) pp. 85--210.

\bibitem{pedregal}
{\sc Pedregal, P.}: {\it Parametrized Measures and Variational Principles.}
Birk\"auser, Basel, 1997.



\bibitem{r}
{\sc Roub\'{\i}\v{c}ek, T.}: {\it Relaxation in Optimization Theory and
Variational Calculus}. W. de Gruyter, Berlin, 1997.

\bibitem{roubicek-kruzik}
{\sc Roub\'{\i}\v{c}ek, T., Kru\v{z}\'{\i}k, M.}: Microstructure evolution model in micromagnetics. {\it Zeit. Angew. Math. Phys.} {\bf 55} (2004), 159--182.

\bibitem{roubicek-kruzik-2}
{\sc Roub\'{\i}\v{c}ek, T., Kru\v{z}\'{\i}k, M.}: Mesoscopical model for ferromagnets with isotropic hardening.  {\it Zeit. Angew. Math. Phys.}  {\bf 56} (2005), 107--135.

\bibitem{schonbek}
{\sc Schonbek, M.E.}: Convergence of solutions to nonlinear dispersive
equations. {\it Comm. in Partial Diff. Equations} {\bf 7} (1982), 959--1000.

\bibitem{silhavy}
{\sc \v{S}ilhav\'{y}, M.}: {\it The Mechanics and Thermodynamics of Continuous Media.} Springer, Berlin, 1997.

\bibitem{silhavy1}
{\sc \v{S}ilhav\'{y}, M.}: Phase transitions with interfacial energy: Interface Null Lagrangians, Polyconvexity, and Existence. In: K.~Hackl (ed.) {\it IUTAM Symposium on Variational Concepts with Applications to the Mechanics of Materials}.
IUTAM Bookseries {\bf 21} (2010), pp.~233--244.
\bibitem{sprenger}
{\sc Sprenger, P.}: {Quasikonvexit\"{a}t am Rande und Null-Lagrange-Funktionen  in der nichtkonvexen Variationsrechnung.} Ph.D. Thesis, Universit\"{a}t Hannover, 1996. 


\bibitem{tartar}
{\sc Tartar, L.}: Compensated compactness and applications to partial
differential equations. In: {\it Nonlinear Analysis and Mechanics} (R.J.Knops,
ed.) Heriott-Watt Symposium IV, Pitman Res. Notes in Math. {\bf 39}, San
Francisco, 1979.

\bibitem{tartar1}
{\sc Tartar, L.}: Mathematical tools for studying oscillations and concentrations: From Young measures to $H$-measures and their variants. In: {\it Multiscale problems in science and technology. Challenges to mathematical analysis and perspectives.} (N.Antoni\v{c} et al. eds.) Proceedings of the conference on multiscale problems in science and technology, held in Dubrovnik, Croatia, September 3-9, 2000.  Springer, Berlin,  2002.

\bibitem{valadier}
{\sc Valadier, M.}: Young measures. In: {Methods of Nonconvex Analysis}
(A.Cellina, ed.) Lecture Notes in Math. {\bf 1446}, Springer, Berlin, 1990,
pp. 152--188.

\bibitem{warga}
{\sc Warga, J.}: {\it Optimal Control of Differential and Functional
Equations.} Academic Press, New York, 1972.

\bibitem{y}
{\sc Young, L.C.}: Generalized curves and the existence of an attained
absolute minimum in the calculus of variations. {\it Comptes Rendus de la
Soci\'et\'e des Sciences et des Lettres de Varsovie}, Classe III {\bf 30}
(1937), 212--234.


}

\end{thebibliography}
\end{document}